%% file: BN2nagdef.tex
\documentclass[10pt]{article}
\usepackage{amssymb,amsmath,stmaryrd}
\usepackage[all]{xy}
\input{engladef.tex}

\def\Dbar{\leavevmode\lower.6ex\hbox to 0pt{\hskip-.23ex
    \accent"16\hss}D}
\def\cfac#1{\ifmmode\setbox7\hbox{$\accent"5E#1$}\else
    \setbox7\hbox{\accent"5E#1}\penalty 10000\relax\fi\raise 1\ht7
    \hbox{\lower1.15ex\hbox to 1\wd7{\hss\accent"13\hss}}\penalty 10000
    \hskip-1\wd7\penalty 10000\box7}
\def\cftil#1{\ifmmode\setbox7\hbox{$\accent"5E#1$}\else
    \setbox7\hbox{\accent"5E#1}\penalty 10000\relax\fi\raise 1\ht7
    \hbox{\lower1.15ex\hbox to 1\wd7{\hss\accent"7E\hss}}\penalty 10000
    \hskip-1\wd7\penalty 10000\box7}
\def\fr{\mathop{\mathrm{Fr}}\nolimits}
\def\sol{\mathop{\mathrm{sol}}\nolimits}
\def\Sol{\mathop{\mathrm{Sol}}\nolimits}
\def\Tr{\mathop{\mathrm{Tr}}\nolimits}
\def\Lti{\overline{L^\times}}
\def\Ltim#1{\overline{L^\times_#1}}

\def\Qtim#1{\overline{\QM^\times_{p,#1}}}
\def\ldots{...}
\title{On modified circular units \\ and annihilation of real classes
}
\author{Jean-Robert Belliard\and Th{\cfac{o}}ng
Nguy{\cftil{e}}n-Quang-{\Dbar}{\cftil{o}}}
\date{january 24 2005}
\begin{document}
\footnotetext{MSC number primary : 11R23, 11R18, 11R20}
\bibliographystyle{alpha}
\maketitle
\setcounter{section}{-1}
\begin{abstract} For an abelian totally real number field $F$ and an
odd prime number $p$ which splits totally in $F$, we present a functorial
approach to special ``$p$-units'' previously built by D. Solomon
using ``wild'' Euler systems. This allows us to prove a conjecture
of Solomon on the annihilation of the $p$-class group of $F$ (in the
particular context here), as well as related annihilation results
and index formulae.
\end{abstract}
\section{Introduction}
Let $F$ be an imaginary abelian field, $G=\Gal(F/\QM)$, and let
$\Cl_F^{\pm}$ be the "plus" and "minus" parts defined by complex
conjugation acting on the class group $Cl_F$ of $F$. The classical
Stickelberger theorem asserts that the Stickelberger ideal, say
$\stc (F)\subset \ZM [G]$, is an annihilator of $\Cl_F^-$. Sinnott
has shown that the index $(\ZM [G]^- : \stc(F)^-)$ is essentially
"equal" to the order $h_F^-$ of $\Cl_F^-$. Using the $p$-adic
point of view, these results fit into the theory of Iwasawa and of
$p$-adic $L$-functions via what is still called the Main
Conjecture, even though it is now a theorem (Mazur-Wiles).\par For
$\Cl_F^+$, our knowledge is more fragmented. For simplicity take
$p\neq 2$. The Main Conjecture relates (in a ``numerical'' way)
the $p$-part $X_F^+$ of $Cl_F^+$ to the $p$-part of
$U_{F^+}/C_{F^+}$, where $U_{F^+}$ (resp. $C_{F^+}$) is the group
of units (resp. circular units) of the maximal real subfield $F^+$
of $F$. In the semi-simple case, it also gives an ideal $MW(F)$ which annihilates a certain
part of $X_F^+$ (the maximal subgroup acting trivially on the
intersection of the $\ZM_p$-cyclotomic extension of $F^+$ with the
$p$-Hilbert class field of $F^+$). In the papers \cite{So92},
\cite{So94} Solomon constructed what could be conside\-red as real
analogues of Gau\ss\  sums and of Stickelberger's element, from
which he conjectured  an annihilation of $X_F^+$. Before sketching
Solomon's construction let us recall some main principles of the
demonstration of Stickelberger's theorem using Gau\ss\ sums : for
all prime ideals $\wp$ not dividing the conductor of $F$, one constructs the
Gau\ss\ sum $g(\wp)$ (which belongs to some cyclotomic extension
of $F$), and if $x$ is a denominator of the Stickelberger element
$\theta_F \in \QM [G]$, then $(g(\wp)^x)=\wp^{x\theta_F}$. We
stress two important ingredients :
\begin{itemize}
\item one can always choose in the class of $\wp$ prime ideals
$\LG$ which split in $F/\QM$.
\item the Gau\ss\ sum is an $\ell$-unit ($\LG \mid \ell$), and
the computation of its $\LG$-adic valuation plays a crucial role.
\end{itemize}
We return to Solomon's construction. For simplicity, we now assume
that $F$ is a (totally) real field, and $p$ an odd prime number
splitting in $F$. From a fixed norm coherent sequence of
cyclotomic numbers (in the extension $F(\mu_{p^\infty})$), Solomon
constructs, using a process which could be seen as a "wild"
variant of Kolyvagin-Rubin-Thaine's method, a special "$(p)$-unit"
$\kappa(F,\gamma)$ (which depends on the choice of a topological
generator of $\Gal (F(\mu_{p^\infty})/ F(\mu_p))$, $\ga$ say). In
fact $\kappa(F,\gamma)\in U^\prime_F\otimes \ZM_p$ where
$U^\prime_F$ is the group of $(p)$-units of $F$. The main result
of \cite{So92} describes the $\PG$-adic valuation of
$\kappa(F,\gamma)$. Let $f$ be the conductor of $F$, put
$\zeta_f=\exp(2i\pi/f)$, and define
$\ep_F=N_{\QM(\zeta_f)/F}(1-\zeta_f)$. Then for $\PG \mid p$ we
have the $p$-adic equivalence :
$$v_\PG (\kappa(F,\gamma))\sim\frac 1 p \log_p (\imath_\PG(\ep_F))$$
where $\imath_\PG$ is the embedding $F\hookrightarrow \QM_p$
defined by $\PG$.
 From this computation and fixing an embedding $\imath \colon
F\hookrightarrow \QM_p$,  Solomon introduces the element
$$\sol_F:=\frac 1 p \sum_{g\in G} \log_p (\imath(\ep_F^g))g^{-1}
\in \ZM_p [G], $$
and conjectures that this element annihilates the $p$-part $X_F$ of the
class group of $F$
(\cite{So92}, conjecture 4.1; actually Solomon states a more general
conjecture assuming
only that $p$ is unramified in $F$).
Here some comments are in order :
\begin{itemize}
\item Solomon's construction is not a functorial one, which means that
he is dea\-ling with elements instead of morphisms (this objection
may be addressed to most "Euler system" constructions); this may
explain the abundance of technical computations required to show
the main result or some (very) special cases of the conjecture.
\item In the semi-simple case (i.e. $p \nmid [F:\QM]$), Solomon's
conjecture is true
but it gives nothing new since one can easily see that
$\sol_F\in MW(F)$ in that case.
\end{itemize}
The object of this paper is to present Solomon's construction in a
more functorial way, making a more wholehearted use of Iwasawa
theory. This functorial approach is not gratuitous,  since on the
one hand it is an alternative to the technical computations in
\cite{So92},\cite{So94}, and on the other hand it allows us to
prove a slightly modified version of Solomon's conjecture and
other related annihilation results
in the most important special case, namely when $p$ splits in
$F$ (see theorem \ref{theo}).\par
We now give an idea of our functorial construction and of
its main derivatives~:\par \noindent Assuming the above mentioned
hypotheses ($F$ is totally real and $p$ splits in $F$), following
$U_F$ and  $C_F$, we introduce $U^\prime_F$, the group of
$(p)$-units of $F$, together with the $p$-adic completions $\U_F,\
\C_F$ and $\U_F^\prime$. The natural idea would be to apply
Iwasawa co-descent to the corresponding $\Lambda$-modules
$\U_\infty,
 \ \C_\infty,\ \U_\infty^\prime$, obtained by taking inverse limits
up the cyclotomic $\ZM_p$-extension of $F$,  but since we are
assuming that $p$ splits in $F$, some natural homomorphisms, for
example $(\C_\infty)_\Gamma \longrightarrow \C_F $, are trivial!
The main point will be to construct other natural (but non
trivial) morphisms, which will allow us to compare the modules at
the level of $F$ with the corresponding modules arising from
Iwasawa co-descent. To fix the ideas let us see how to compare
$\C_F$ with $(\C_\infty)_\Gamma$. We construct successively two
homomorphisms :\par
\noindent $1)$ The first one, whose
construction essentially uses an analog of "Hilbert's theo\-rem 90
in Iwasawa theory", allows us to map $(\C_\infty)_\Gamma$ to
$\U^\prime_F$ after ``dividing by $(\ga -1)$'' (see theorem \ref{iso}).
One may even go to the quotient by $\U_F$, the
natural projection being injective when restricted to the image of
$(\C_\infty)_\Gamma$. So we have canonical maps~:\par $\xymatrix {
(\C_\infty)_\Gamma \ar[r]^-\alpha& \U^\prime_F}$ and $\xymatrix {
(\C_\infty)_\Gamma \ar[r]^-\beta& \U^\prime_F/\U_F}$. Then the
image under $\be$ of an ad hoc element in $\C_\infty$ is nothing
else but Solomon's element $\kappa(F,\ga)$ and functorial
computations with $\be$ allow us to solve all index questions
arising in \cite{So92},\cite{So94}.\par \noindent $2)$ The second
homomorphism constructs a "bridge" linking Sinnott's exact
sequence with the class-field exact sequence. More precisely,
$\pucod_F:=\al((\C_\infty)_\Ga)$ is considered as a submodule of
$\U_F^\prime/ \U_F$ as in $1)$, while $\C_F$ is considered as a
submodule of $\U^\prime_F/\Uha_F^\prime$, where $\Uha_F^\prime$ is
some kernel consisting of local cyclotomic norms, usually called
the Sinnott-Gross kernel. A priori there doesn't exist any natural
map from $\U_F^\prime/\Uha_F^\prime$ to   $\U_F^\prime/\U_F$,
since $\U_F\bigcap\Uha_F^\prime = \{1\}$,
but when restricting to $\C_F$ and using the map $\be$, one is able to
construct a map $\varphi \colon \C_F \longrightarrow \ _{\Tr}
(\pucod_F)$ (the kernel of the algebraic trace of $\ZM_p [G]$
acting on $\pucod_F$), which is injective and whose cokernel may
be computed (see theorem \ref{epimon}) .\par
Once the morphisms $\alpha, \beta, \varphi$ have
been constructed one is able to compute the Fitting ideals in
$\ZM_p [G]$ of various codescent modules in terms of Solomon's
elements. This yields annihilation results for various classes,
e.g. ideal classes, but also quotients of global units or
semi-local units modulo circular units (see theorem \ref{fit}).
Note that in the non
semi-simple case, these differ from the ones obtained by means of
usual (tame) Euler systems.

\par
In recent developments (\cite{BB01}, \cite{BG00}, \cite{RW03})
around the ``lifted  root number conjecture'' and its ``twisted''
versions (or ``equivariant  Tamagawa number'' conjecture for Tate motives
over $\QM$ in the language of Burns and Flach), Solomon's elements appear
to play a  role. In some special cases, the LRN conjecture is  equivalent to
the existence of $S$-units satisfying a variety of  explicit
conditions; in particular, it implies a refinement of Solomon's main result
(op. cit.) on the $\PG$-adic valuation of Solomon's element  (\cite{BB01},
3.2 and 4.14). Moreover, the proof of the ETN conjecture for abelian
fields proposed recently by \cite{BG00} uses Iwasawa co\--descent in its final
step and, in the splitting case, Solomon's elements appear again
(\cite{BG00}, 8.7, 8.11, 8.12, 8.13), together
with the usual heavy calculations (see also \cite{RW03}). Actually,
all  the technicalties seem to be related to the phenomenon of ``trivial
zeroes'' of $p$-adic $L$-functions (with an intended vagueness, this means~:
zeroes of local Euler factors at $p$). This renders our present approach all
the more interesting, since it appears to give a functorial (as opposed to
technical) process to ``bypass trivial zeroes''. We hope to come back
to this topic in a subsequent work.

\section{A few exact sequences}\label{SEs}
In this section (which can be skipped at first reading), we collect
for the convenience of the reader a few ``well known'' exact sequences
which come from class-field theory and Iwasawa theory, and will be used freely
in the rest of the paper. Let $F$ be a number field, $p$ an odd prime,
$S=S(F)$ the set of primes in $F$ which divide $p$. We 'll adopt once and for
all the following notations :
\begin{itemize}
\item[] $U_F$ (resp. $U^\prime_F$) $=$ the group of units (resp. $S$-units)
of $F$.
\item[] $X_F$ (resp. $X_F^\prime$) $=$ the $p$-group of ideal classes (resp.
$S$-ideal classes) of $F$.
\item[] $\mathfrak {X}_F \ =$ the Galois group over $F$ of the maximal
 $S$-ramified (i.e. unramified outside $S$) abelian pro-$p$-extension of
$F$.
\item[] For any abelian group, we'll denote by $\bar A= \underset {\underset m
 \leftarrow} \lim\ A/A^{p^m}$ the $p$-completion of $A$. If $A$ is of finite
type, then $\bar A = A \otimes \ZM_p$.
\end{itemize}
\subsection{Exact sequences from class-field theory}\label{SEcdc}
By class-field theory, $X_F$ (resp. $X^\prime_F$) is canonically isomorphic
to the Galois group over $F$ of the maximal unramified (resp. unramified and
$S$-decomposed) abelian extension $H_F$ (resp. $H_F^\prime$) of $F$. We have
two exact sequences :
\begin{itemize}
\item[--] one relative to inertia :\par
\noindent \hfill \xymatrix {\U_F
\ar[r]^-{\mathrm {diag}} & \mathcal U := \bigoplus_{v\in S} U^1_v
\ar[r]^-{\mathrm {Artin}} &
\mathfrak  {X}_F\ar[r]& X_F\ar[r]&0}
\hfill $(1)$\par
Here $U^1_v=\U_{F_v}$ is the group of principal local units at $v$,
and ``diag'' is induced by the diagonal map. Assuming Leopoldt's conjecture
for $F$ and $p$ (which holds e.g. when $F$ is abelian), it is well known
that ``diag'' is injective.
\item[--] one relative to decomposition :\par
\noindent \hfill \xymatrix {\U^\prime_F
\ar[r]^-{\mathrm {diag}} & \mathcal F := \bigoplus_{v\in S} \overline
{F^\times_v}
\ar[r]^-{\mathrm {Artin}} &
\mathfrak  {X}_F\ar[r]& X^\prime_F\ar[r]&0}
\hfill $(2)$\par
Here again, the map ``diag'' is injective modulo Leopoldt's conjecture.
By putting together $(1)$ and $(2)$ (or by direct computation), we get
a third exact sequence :\par
\noindent \hfill \xymatrix {0\ar[r] & \U^\prime_F/\U_F
\ar[r]^-\mu & \mathcal F /\mathcal U
\ar[r]^-{\mathrm {Artin}} &
X_F\ar[r]& X^\prime_F\ar[r]&0}
\hfill $(3)$
\end{itemize}
where $\mu$ is induced by the valuation maps at all $v\in S$. By local class
field theory $\overline {F^\times_v}/ U^1_v$ is canonically isomorphic to
the Galois group over $F_v$ of the unramified $\ZM_p$-extension of $F_v$.
If $F$ is Galois over $\QM$, then $\mathcal F /
\mathcal U \simeq \ZM_p [S]$ as Galois modules, and we can rewrite $(3)$ as :
\par
\noindent \hfill \xymatrix {0\ar[r]&\U^\prime_F/\U_F
\ar[r]^-\mu & \ZM_p [S]
\ar[r]^-{\mathrm {Artin}} &
X_F\ar[r]& X^\prime_F\ar[r]&0}
\hfill ($3\ \mu$)
\subsection{Sinnott's exact sequence}\label{Si}
For any $v\in S$, let $\widehat {F}_v^\times$ be the subgroup of
 $\overline {F^\times_v}$ which corresponds by class-field theory to the
cyclotomic $\ZM_p$-extension $F_{v,\infty}$ of $F_v$ i.e.
$\overline {F^\times_v}/\widehat {F}_v^\times\simeq \Gal(F_{v,\infty}/F_v)$
and $\widehat F_v^\times$ is usually called the group of {\it universal
(cyclotomic) norms} of $F_{v,\infty}/F_v$.
Let $\de \colon \U^\prime_F \longrightarrow \oplus_{v\in S}
\overline {F^\times_v}/\widehat F_v^\times$ be the homomorphism induced
by the diagonal map. Its kernel, $\widehat U^\prime_F$ say, usually called
the {\it Gross kernel} (\cite{FGS}; see also \cite{Ku72}), consists of
elements of $\U_F^\prime$ which are everywhere universal cyclotomic norms.
Its cokernel is described by the Sinnott exact sequence (see the appendix to
\cite{FGS}), for which we need additional notations.\par
Let $F_\infty=\bigcup_{n\geq 0} F_n$ be the cyclotomic $\ZM_p$-extension of
$F$, $L_\infty^\prime$ the maximal unramified abelian
pro-$p$-extension of $F_\infty$
which is totally decomposed at all places dividing $p$ (hence at all places),
$L^\prime_0$ the maximal abelian extension of $F$ contained in
$L^\prime_\infty$. Then Sinnott's exact sequence reads :\par
\noindent \hfill \xymatrix {0 \ar[r] & \U^\prime_F/\widehat {U}^\prime_F
\ar[r]^-{\de} &  \bigoplus_{v\in S}
\overline {F^\times_v}/\widehat F_v^\times
\ar[r]^-{\mathrm {Artin}} &
\Gal(L^\prime_0/F)\ar[r]& X_F^\prime\ar[r]&0}
\hfill $(4)$\par
By class-field theory, $\Gal(L^\prime_\infty/F_\infty)\simeq X^\prime_\infty
:= \underset \leftarrow \lim \
X^\prime_{F_n}$ and $\Gal(L^\prime_0/F_\infty)\simeq
(X_\infty^\prime)_\Ga$ (the co-invariants of $X^\prime_\infty$ by $\Ga =
\Gal(F_\infty/F)$). By the product formula, the image of $\de$ is contained
readily in $\widetilde \oplus_{v \in S}
\overline {F^\times_v}/\widehat F_v^\times$ (:= the kernel of the map ``sum
of components'') and $(4)$ can be rewritten as :\par
\noindent \hfill \xymatrix {0 \ar[r] & \U^\prime_F/\widehat {U}^\prime_F
\ar[r]^-{\de} &\widetilde  \bigoplus_{v\in S}
\overline {F^\times_v}/\widehat F_v^\times
\ar[r]^-{\mathrm {Artin}} &
(X^\prime_\infty)_\Ga\ar[r]& X_F^\prime}
\hfill $(5)$\par
The image of the natural map $(X^\prime_\infty)_\Ga\longrightarrow X_F^\prime $
is nothing but $\Gal(H^\prime_F/H^\prime_F\cap F_\infty)$. Gross' conjecture
asserts that $(X^\prime_\infty)_\Ga$ is finite. It holds (which is the case if
$F$ is abelian) if and only if $\Uha^\prime_F$ has $\ZM_p$-rank equal to
$r_1+r_2$. Let us look more
closely at the Galois structure of $\oplus_{v\in S}
\overline {F^\times_v}/\widehat F_v^\times$. For simplicity,
suppose that the completions of $F$ are
linearly disjoint from those of $\QM_\infty$ (this property should
 be called ``local linear disjointness'')
: this happens e.g. if $p$
is totally split in $F$. Then the local norm map gives an isomorphism
$\overline {F^\times_v}/\widehat {F_v^\times}\simeq \overline {\QM_p^\times}
/\widehat \QM_p^\times$; besides $\overline {\QM_p^\times}
/\widehat \QM_p^\times\simeq 1+p\ZM_p$ is isomorphic to $p\ZM_p$ via the
 $p$-adic logarithm, hence $\oplus_{v\in S}
\overline {F^\times_v}/\widehat F_v^\times\simeq p \ZM_p[S]$ in this situation.
Let us denote by $I(S)$ the kernel of the sum in $\ZM_p[S]$. Assuming that
$p$ is totally split in $F$, we can rewrite $(5)$ as :\par
\noindent \hfill \xymatrix {0 \ar[r] & \U^\prime_F/\widehat {U}^\prime_F
\ar[r]^-{\la} & pI(S)
\ar[r] &
(X^\prime_\infty)_\Ga\ar[r]& X_F^\prime\ar[r] & 0}
\hfill ($5\ \la$)\par
On comparing the exact sequences ($3\ \mu$) and ($5\ \la$), it is tempting
to try and ``draw a bridge'' between them, but unfortunately, there seems to be
no natural link between $\U^\prime_F/\widehat U^\prime_F$ and
$\U^\prime_F/\U_F$ because $\widehat U^\prime_F$ and $\U_F$ are
``independent'' in the following sense :
\begin{lem}\label{cdcl} Suppose that $F$ satisfies Leopoldt's conjecture and
$p$ is totally split in $F$. Then
$$ \Uha_F^\prime\cap \U_F =\{1\}$$
\end{lem}
\dem Assuming Leopoldt's conjecture, we can embed $\U_F$ and
$\widehat U^\prime_F$ in $\mathcal F = \oplus_{v\in S} \overline {F_v^\times}$.
If $p$ is totally split, then for any $v\in S$,
$\widehat{F}^\times_v=\widehat {\QM}^\times_p=p^{\ZM_p}$
and $U^1_v=1+p \ZM_p$, hence an element
$u\in \U_F$ belongs to $\widehat U^\prime_F$ if and only if $u=1$.
\par\qed\par
Nevertheless, one of our main results will be the construction of a natural
map ($\S$ \ref{bridge}) between two appropriate submodules of
$\U^\prime_F/\widehat U^\prime_F$ and
$\U^\prime_F/\U_F$.
\section{Hilbert's theorem 90 in Iwasawa theory}\label{hi90}

\subsection{Some freeness results}\label{freres}

We keep the notations of $\S$ \ref{SEs}. Let $\La = \ZM_p[[\Ga]]$ be the
Iwasawa algebra.
If $ M_\infty=\limpro M_n $ is a $\Lambda$-module,
$\alpha_\infty=(\alpha_0,\ldots,\alpha_n, \ldots)$ will denote
a typical element of $M_\infty$, the index zero
referring to $F=F_0$. The canonical map $(M_\infty)_\Gamma
\longrightarrow M_0$ sends $\alpha_\infty+I_\Gamma M_\infty$ (where
$I_\Ga$ is the augmentation ideal) to the component $\alpha_0$
of $\al_\infty$. Let us denote by $M_\infty^{(0)}$ the submodule of
$M_\infty$ consisting of the elements
$\alpha_\infty$ such that $\alpha_0=0$ (in
additive notation). We want to study the relations between the $\La$-modules
$\U_\infty:=\limpro \U_n$ and $ \U^\prime_\infty:=\limpro \U^\prime_n$.
Obviously, they have the same $\La$-torsion, which is $\ZM_p(1)$ or $(0)$
according as $F$ contains or not a primitive $p^{\mathrm {th}}$ root of
unity (see e.g. \cite{Ku72} or \cite{Wi85}).
Let $\fr_\La (\U_\infty) := \U_\infty/\tor_\La(\U_\infty)$
and $\fr_\La (\U^\prime_\infty) := \U^\prime_\infty/
\tor_\La(\U^\prime_\infty)$.
\begin{prop}\label{Ufree}
 Let $F$ be any number field, $[F:\QM]=r_1+2r_2$.
Then the $\La$-modules $\fr_\La (\U_\infty)$ and
$\fr_\La (\U^\prime_\infty)$ are free, with $\La$-rank
equal to $r_1+r_2$.
\end{prop}
\dem The assertion concerning $\fr_\La (\U^\prime_\infty)$
is a theorem of Kuz$'$min (\cite{Ku72}, 7.2). As for $\fr_\La
 (\U_\infty)$, it is enough to notice that the quotient
$\fr_\La (\U^\prime_\infty)/\fr_\La (\U_\infty) \simeq
\U^\prime_\infty/\U_\infty$ has no $\ZM_p$-torsion, hence has no
non-trivial finite submodule. This is equivalent to saying that
$\fr_\La (\U^\prime_\infty)/\fr_\La (\U_\infty)$ is of projective
dimension at most $1$ over $\La$ (see e.g. \cite{Ng84}).
Since $\La$ is local, the $\La$-freeness of $\fr_\La
(\U_\infty)$ follows by Schanuel's lemma.\par\qed

\subsection{Cyclotomic submodules}\label{cycfree?}

If the field $F$ is abelian over $\QM$, we have at our disposal
the group $C_F$ (resp. $C^\prime_F$) of {\it Sinnott's circular
units} (resp. {\it circular $S$-units}), which is defined
as being the intersection of $U_F$ (resp. $U^\prime_F$)
with the group of {\it circular numbers} of $F$
(\cite{Si78}, \cite{Si80}; see also \S \ref{cycnum} below). Sinnott's index
formula states that $(U_F:C_F)=c_F h_F^+$, where $h^+_F$ is
the class number of the maximal real subfield of $F$, and $c_F$ is
a rational constant whose explicit definition does not involve the class
group. Going up the cyclotomic $\ZM_p$-extension $F_\infty=\bigcup_n F_n$,
it is known that the cons\-tants $c_{F_n}$ remain bounded.
It is a classical result (see e.g. \cite{Grei92}) that $\C_\infty$ and
$\U_\infty$ have the same $\La$-rank. It is an obvious consequence of
the definition that $\C_\infty$ and $\C_\infty^\prime$ have the same
$\La$-rank.
If the base field $F$ is a cyclotomic field, it is also known,
by results on distributions ``\`a la Kubert-Lang'', that
$\fr_\La (\C_\infty)$ and  $\fr_\La (\C^\prime_\infty)$
are $\La$-free (see e.g. \cite{Ku96}). Life would be too easy if
such results could be extended to any abelian field.
By proposition \ref{Ufree} and the exact sequences
\xymatrix{0\ar[r] & \fr_\La(\C_\infty) \ar[r] &
\fr_\La(\U_\infty) \ar[r] & \U_\infty/\C_\infty \ar[r] & 0}, and
\xymatrix{0\ar[r] & \fr_\La(\C^\prime_\infty) \ar[r] &
\fr_\La(\U^\prime_\infty) \ar[r] & \U^\prime_\infty/
\C^\prime_\infty \ar[r] & 0},
we see that the $\La$-freeness of $\fr_\La(\C_\infty)$
(resp $\fr_\La(\C^\prime_\infty)$) is
equivalent to the triviality of the maximal finite submodule
of $\U_\infty/\C_\infty$ (resp. $\U^\prime_\infty/\C^\prime_\infty$).
Let us call $MF_\infty$ (resp. $MF^\prime_\infty$) these maximal
finite submodules. In order to get hold of $MF_\infty$ and
$MF_\infty^\prime$, let us define~:
\begin{defi}
 The {\rm Ku{\v{c}}era-Nekov{\'a}{\v{r}} kernel}
is defined by
$$KN_{F,n}:=\bigcup_{m\geq n} \Ker i_{n,m}$$
where for $m\geq n \in \NM$,
$i_{n,m} \colon \U_n/\C_n \longrightarrow \U_m/\C_m$ is the natural map.
\end{defi} The orders of the kernels
$KN_{F,n}$ are bounded independently of $n$ :
this is an immediate consequence of the main result of \cite{GK89}
and the one of \cite{KN95}.
Let us denote by $KN_F$ the projective limit (relatively to norm maps
and $n$) $KN_F:=\limpro KN_{F,n}$. This kernel $KN_F$ is
the obstruction to the $\La$-freeness of $\fr(\C_\infty)$ and
$\fr(\C^\prime_\infty)$,
in the following sense~:
\begin{prop}\label{KNF}\ \par
\begin{itemize}
\item[{\rm $(i)$}] $(\U_\infty/\C_\infty)^\Ga$ is finite.
\item[{\rm $(ii)$}] We have equalities
$KN_F= MF_\infty= MF_\infty^\prime$, and canonical isomorphisms
$(KN_F)^\Ga \simeq \tor_{\ZM_p} (\fr(\C_\infty))_\Ga\simeq \tor_{\ZM_p}
(\fr(\C^\prime_\infty))_\Ga$. In particular, $\fr_\La (\C_\infty)$
(resp. $\fr_\La (\C^\prime_\infty)$) is free if and only if
$KN_F=0$.\end{itemize}
\end{prop}
\dem $(i)$
The Main Conjecture (or Mazur-Wiles' theorem) applied to the
ma\-xi\-mal real subfield $F^+$ of $F$
implies that the $\La $-torsion modules
$\U^+_\infty/\C^+_\infty=\U_\infty/\C_\infty$ and
$X_\infty^+$ have the same characteristic series
(where the
$ ^+$ denotes the objects related to $F^+$, i.e. the submodule
on which complex
 conjugation acts trivially because $F$ is abelian and
 $p\neq 2$). Because
 Leopoldt's conjecture holds all along the cyclotomic
tower, it is well known that $X_\infty^+$ and $(X_\infty^\prime)^+$ are
pseudo-isomorphic (see e.g. \cite{Wi85}). By Gross~'~conjecture (which
holds because $F$ is abelian),
 $((X_\infty^\prime)^+)^\Ga$
is finite, and so is
$(\U_\infty/\C_\infty)^\Ga$.\par
$(ii)$ We have an exact sequence :\par
\noindent\centerline{\xymatrix{0\ar[r] & \U_\infty/\C_\infty \ar[r] &
\U^\prime_\infty /\C^\prime_\infty\ar[r] & \U^\prime_\infty/(\U_\infty
+ \C_\infty^\prime)\ar[r] & 0}}\par
\noindent It gives an inclusion $MF_\infty\subset MF_\infty^\prime$.
To prove the inverse inclusion it is enough to show that $
\U^\prime_\infty/(\U_\infty
+ \C_\infty^\prime)$ is without $\ZM_p$-torsion. Let $\SC$ be the (finite) set
of places above $p$ of $F_\infty$. Clearly $\U^\prime_\infty/\U_\infty$
is isomorphic to a submodule of finite index of
$\ZM_p[\SC]$ with the natural action of $\Gal(F_\infty/\QM)$
on both sides. Call $M$ this submodule.
It follows that $\U^\prime_\infty/(\U_\infty+\C^\prime_\infty)$
is isomorphic to $M/\mu(\C^\prime_\infty)$, where $\mu$
is obtained from the valuations at finite levels, i.e. is the
limit of the homomorphisms $\mu$ of the exact sequences $(3)$.

Let $\BM_\infty$ be the cyclotomic $\ZM_p$-extension of $\QM$, and $\BM_n$ its
$n^{\text th}$ step.
Let $\ep_{\QM,\infty}$ denotes the norm coherent sequence of numbers
$\ep_{\QM,\infty}=(p,N_{\QM(\ze_{p^2})/\BM_1}(1-\ze_{p^2}), \dots,
N_{\QM(\ze_{p^{n+1}})/\BM_n} (1-\ze_{p^{n+1}}),\dots )$. Since
$F_\infty$ is abelian over $\QM$ it is at most tamely ramified over
$\BM_\infty$, with ramification index $e$ say. It follows that
$\mu(\ep_{\QM,\infty})=e  \sum_{v\in \SC} v$, with $e\in \NM$ dividing
$(p-1)$ and therefore a unit in $\ZM_p$. Each generators of
$\C_\infty$ either is a  unit or behaves like $\ep_{\QM,\infty}$.
This shows that $\mu(\C^\prime_\infty)=\ZM_p \sum_{v\in \SC} v$.
Consequently we have $\U^\prime_\infty/(\U_\infty+\C^\prime_\infty)
\simeq M/(\ZM_p \sum_{v\in \SC} v)\subset
\ZM_p[\SC]/(\ZM_p \sum_{v\in \SC} v)$ and this proves
the second equality of proposition \ref{KNF} $(ii)$. \par
For the first equality, there exist an $n_0$ such that
for all $n\geq n_0$, $\ga^{p^n}$ acts trivially on $MF_\infty$.
Putting $\omega_n = \ga^{p^n}-1$ as usual, and taking $n\geq n_0$ and $m-n$
large enough in the commutative triangle~:\par
\noindent\centerline{\xymatrix{ \U_m/\C_m \ar[rr]^-{\omega_m/\omega_{n}}
\ar[rd]_-{N_{m,n}} & & \U_m/\C_m \\ & \U_{n}/\C_{n} \ar[ru]_-{i_{n,m}} }}
where $N_{m,n}$ is the obvious norm map, we see immediately that
$MF_\infty \subset KN_F$, and the finiteness of $KN_F$ shows the equality.
 Moreover, the exact sequence of $\La$-modules  \par
\centerline{\xymatrix{0\ar[r] & \fr_\La(\C_\infty) \ar[r] &
\fr_\La(\U_\infty) \ar[r] & \U_\infty/\C_\infty \ar[r] & 0}}
\par
\noindent  gives
by descent an exact sequence of $\ZM_p$-modules~:
\par
\noindent\centerline{\xymatrix@=22pt{0\ar[r] & (\U_\infty/\C_\infty)^\Ga
\ar[r] & \fr_\La(\C_\infty)_\Ga \ar[r] &
\fr_\La(\U_\infty)_\Ga\ar[r] & (\U_\infty/\C_\infty)_\Ga \ar[r] & 0 }}
\par
\noindent By $(i)$ $(\U_\infty/\C_\infty)^\Ga$
is finite, hence the equality $(\U_\infty/\C_\infty)^\Ga=(MF_\infty)^\Ga$.
It follows then from \ref{Ufree} that $(MF_\infty)^\Ga=\tor_{\ZM_p}
(\fr_\La (\C_\infty)_\Ga)$. The remaining assertions in \ref{KNF}
are straightforward.\par\qed\par
In what follows the finite order $\# (KN_F)^\Ga$ will play
for descent modules a role analogous to the one played
by Sinnott's constant $c_F$ at finite level. Let us put the~:
\begin{nota}\label{MFord} $\kappa_F:=\# (KN_F)^\Ga$\end{nota}
\begin{itemize}
\item[] {\it Remark~:}  We take this opportunity to correct a few mistakes
in \S $3$ and \S $4$ of \cite{BN1} (due to the eventual non triviality of
the constant $\kappa_F$). In short~: every statement concerning only
pseudo-isomorphisms or characteristic series remains true; every index
formula should be corrected if necessary by a factor involving $\kappa_F$;
in every monomorphism or isomorphism statement related to $(\C_\infty)_\Ga$,
this module should be replaced by its image in $\U_F$. Another, quicker
solution would be to restrict generality and
suppose everywhere that $\kappa_F=1$.
\end{itemize}
\noindent In many situations (e.g.when $p\nmid [F:\QM]$, or
 when $F$ is a cyclotomic
field), it is known that $KN_F=0$, but note
that this is indirect evidence, coming from the $\La$-freeness of
$\fr_\La (\C_\infty)$ (\cite{Ku96}, \cite{JNT1} ...), not
from the definition of the obstruction kernels. Let us describe  examples
of non trivial $KN_F$ (again by
indirect evidence, this time finding a $(\C_\infty)_\Ga$ containing
non trivial $\ZM_p$-torsion). Note that another example for $p=3$ has
just been announced by R. Ku\v cera (\cite{KuJA}).
 Such examples may be considered as exceptional~: let us recall that
we have to avoid {\it Hypoth\`ese B} of \cite{JNT1}, which
holds true in most cases. In order to prove the non-triviality
of the constant $\kappa_F$ we need to
 add some very peculiar decomposition
hypotheses such as those used in \cite{Grei93}. Then
the Galois module structure of circular units is less difficult
to control. This justifies the terminolo\-gy
{\it ``g\"unstige $(p+1)$-tuple''} used in \cite{Grei93}.
We state these conditions :

\begin{itemize}
\item[1--] the conductor $f$ of $F$ is of the form
$f=\prod_{i=1}^{i=p+1} l_i$, where $l_i$ are prime
numbers, $l_i\equiv 1 [p]$, and for all $j\neq i$, there exist
some $x_{i,j}$ such that  $l_i\equiv x_{i,j}^p [l_j]$.
\item[2--] $G:=\Gal(F/\QM) \simeq (\ZM/p\ZM)^2$
\item[3--] All the $(p+1)$ subfields of absolute degree $p$ of $F$
($F^1,F^2,\ldots,F^{p+1}$ say) have conductors
$\cond(F^j)=\prod_{i=1,i\neq j}^{i=p+1} l_i$
\end{itemize}
\begin{itemize}
\item [{}]
\begin{itemize}
\item [{}]{\it Remark~:} Following Greither, and using \v Cebotarev density
theorem, it is not difficult to prove  that there exist infinitely
many $(p+1)$-uples of primes $l_i$ such that $p$ and any subfield
of the cyclotomic field $\QM(\zeta_{\prod l_i})$ satisfy condition 1.
Such $(p+1)$-uples were called {\it g\"unstige $(p+1)$-tuple} in \cite{Grei93}.
We can then deduce that for each $p$ there exist (infinitely) many fields
$F$ such that $F$ and $p$ satisfy condition 1--, 2-- and 3--.
\end{itemize}
\end{itemize}
\begin{prop} If the fixed field $F$ together with the fixed prime
$p\neq 2$ satisfy condition 1--, 2-- and 3--, then
$\fr(\C_\infty)$ is not $\La$-free.
\end{prop}
\dem Details are given in \cite{pmb}.\par\qed

\subsection{Dividing by $T$}\label{h90map}

Going back to the general (not necessarily abelian) situation, let us record
a useful result of co-descent, due to Kuz$'$min~:
\begin{prop}[\cite{Ku72}, 7.3]\label{codes}
The natural map $(\U^\prime_\infty)_\Ga
\longrightarrow \U^\prime_F$ is injective.
\end{prop}
\begin{itemize} \item [{}]
{\it Remark~:} For simple reasons of $\ZM_p$-ranks,
the analogous property for $\U_\infty$ does not hold. In other
words the natural map $(\U_\infty)_\Ga \longrightarrow (\U^\prime_\infty)_\Ga
$ induced by the inclusion $\U_\infty\subset \U^\prime_\infty$, is not
injective. \end{itemize}
We aim to replace this last inclusion by a more elaborate
injection. To this end, we first prove an analog of Hilbert's theorem 90
in Iwasawa theory :
\begin{theo}\label{iso} Recall that $(\overline U_\infty^\prime)^{(0)}$
is defined as the kernel of the naive descent map $\overline
U_\infty^\prime \longrightarrow \overline U_0^\prime$.
Let us fix a topological generator $\ga$ of
$\Ga=\Gal(F_\infty/F)$. Then multiplication by $(\ga -1)$ gives an
isomorphism of $\Lambda$-modules~:

\noindent \begin{center} $ \overline U_\infty^\prime
\overset \sim \longrightarrow (\overline U_\infty^\prime)^{(0)}$. \end{center}
\end{theo}

\dem
\noindent Since $\Ga$ acts trivially on $\U_0^\prime$, the image
$(\ga-1) \U^\prime_\infty$ is obviously contained in
$({\overline U_\infty^\prime})^{(0)}$. Moreover, the kernel of $(\ga-1)$
is $\La$-torsion, but the $\La$-torsion submodule of $\U^\prime_\infty$ is
$(0)$ or $\ZM_p(1)$, hence $(\ga-1)$ is injective on $\U^\prime_\infty$.
It remains to show that $(\ga-1) \U^\prime_\infty =
({\overline U_\infty^\prime})^{(0)}$.
To this end, consider the composite of natural injective maps :
\par
\noindent\centerline{\xymatrix{
({\overline U_\infty^\prime})^{(0)}/(\ga-1) \U^\prime_\infty
\ar@{^{(}->}[r] & \U^\prime_\infty/(\ga-1) \U^\prime_\infty
=(\U^\prime_\infty)_\Ga \ar@{^{(}->}[r]& \U^\prime_F}}
(the injectivity on the left is obvious; on the right it follows
from \ref{codes}). By definition (see the beginning of this section),
this composite map is null, which shows what we want.\par\qed
\begin{cor}\label{dtr} Suppose that $F$ satisfies Leopoldt's conjecture
and $p$ is totally split in $F$. Then multiplication by $(\ga-1)$
gives an isomorphism of $\La$-modules $\U_\infty^\prime \overset \sim
\longrightarrow \U_\infty$
\end{cor}
\dem Using \ref{iso}, it remains only to show the equality
$(\U^\prime_\infty)^{(0)}= \U_\infty$. Let us show successively the mutual
inclusions~:
\par
\noindent $(i)$
If $\alpha_\infty\in
(\U^\prime_\infty)^{(0)}$, then by \ref{iso},
$\alpha_\infty$ may be written as
$\al_\infty=(\gamma-1) \beta_\infty$, with
 $\beta_\infty \in \U^\prime_\infty$. Since $\ga$ acts trivially on the primes
of $F$ above $p$, it is obvious that $(\ga-1)\be_\infty \in \U_\infty$.
Hence $(\U^\prime_\infty)^{(0)}\subset \U_\infty$.
\par
\noindent $(ii)$ Let $u_\infty\in \U_\infty$. For all primes $v$ above
$p$ and all integers $n$, we have $u_0=N_{F_n/F}(u_n)=N_{F_{v,n}/F_v} (u_n)$,
i.e. $u_0$ is actually an element of $\widehat {U}^\prime_F$. But
$\U_F \bigcap \widehat {U}^\prime_F = \{ 1\}$ by \ref{cdcl}. Hence $\U_\infty
\subset (\U^\prime_\infty)^{(0)}$.\par\qed\par
We are now in a position to construct the map we were looking for.
\begin{defi}\label{maps}\ \par
\begin{itemize}
\item[{\rm $(i)$}] Suppose that $F$ satisfies Leopoldt's conjecture and $p$ is
totally split in $F$. Define
\par
\centerline {\xymatrix {\Th\colon \U_\infty\ar[r]^-{\sim}&
\U^\prime_\infty \ar[r]^-{\text {nat.} } &
\left (\U^\prime_\infty \right )_\Gamma \ar@{^{(}->}[r] &\U^\prime_F .}}
\par
The first isomorphism is the inverse of the one in \ref{dtr};
the last monomorphism is the one in \ref{codes}. The definition of $\Th$
depends on the choice of $\ga$.
\item[{\rm $(ii)$}] Suppose that  $F$ is abelian over $\QM$.
Then $F_\infty/\QM$ is also
abelian, and each $F_n$ contains {\rm circular units} in the sense of
Sinnott \cite{Si80}.
Let $C_n$ be the group of circular units of $F_n$, $\C_n=C_n\otimes \ZM_p$,
$\C_\infty = \limpro \C_n$. If $p$ is totally split in $F$, define
$\Th_C \colon \C_\infty \longrightarrow \U^\prime_F$ as the restriction
of $\Th$ to $\C_\infty$.
\end{itemize}
\end{defi}
The homomorphism $\Th_C$ gives us Solomon's ``$S$-units'' at a lower
cost than in \cite{So92}. In fact the choice of any interesting element
in $\C_\infty$ would yield via $\Th_C$ an interesting element
in $\U^\prime_F$. Let us fix once and for all the notation
$(\ze_n)_{n\in \NM}$ for a system of primitive $n^{\text{th}}$ roots of
unity such that $\ze_{mn}^m=\ze_n$. Let $f$ be the conductor of $F$.
Consider the norm coherent sequence of circular units
$\ep_{F,\infty}:=(1, N_{\QM(\zeta_{p^2 f})/F_1} (1-\ze_{p^2 f}), \cdot \cdot
\cdot,N_{\QM(\zeta_{p^{n+1} f})/F_n} (1-\ze_{p^{n+1} f}), \cdot\cdot\cdot
)_{n\in \NM}$. It is an easy (but fastidious) consequence of the definition
that $\Th_C(\ep_{F,\infty})$ is indeed equal to Solomon's ``$p$-unit''
$\kappa(F,c)$, as constructed in \cite{So92}. We skip the calculations
since they give no further information. Actually we don't need to show
the equality between $\kappa (F,c)$ and $\Th_C(\ep_{F,\infty})$ because in
the sequel, all properties of $\Th_C(\ep_{F,\infty})$ which will be
used (e.g. \ref{so3.1}) will be reproved.
\section{Modified circular $S$-units and regulators}\label{indreg}
\subsection{Modified circular $S$-units}
From now on, unless otherwise stated, $F$ will be an
{\sl abelian number field, and $p$ will be totally split in $F$}.
We aim to study the maps $\Th$ and $\Th_C$.
\begin{prop}\label{ThU} $\Ker \Th=(\ga-1) \U_\infty$ and  $\im \Th \simeq
(\U_\infty^\prime)_\Ga$. In particular $\im \Th$ has $\ZM_p$-rank $(r_1+r_2)$.
\end{prop}
\dem The isomorphism $\im \Th \simeq (\U^\prime_\infty)_\Ga$ follows
obviously from the definition of $\Th$ and \ref{codes}. The value
of the $\ZM_p$-rank comes from \ref{Ufree}. To compute $\Ker \Th$, we can
decompose $\Th$ as $\U_\infty \overset {\mathrm{nat}} \twoheadrightarrow
(\U_\infty)_\Ga \overset {\overset {\ref{dtr}} \sim} \rightarrow
(\U^\prime_\infty)_\Ga \overset {\ref{codes}} \hookrightarrow \U^\prime_F$.
It is then obvious that $\Ker \Th = (\ga-1) \U_\infty$. \par\qed\par
It is naturally more difficult to get hold of $\Th_C$, which contains
more arithme\-tical information.
\begin{prop}\label{ThC}
$\Ker \Th_C= ((\ga-1)\U_\infty) \cap \C_\infty$
and (abusing
notation)

\noindent $\im
\Th_C\simeq (\C_\infty)_\Ga / (KN_F)^\Ga$.
In particular
$\im \Th_C$ has $\ZM_p$-rank $(r_1+r_2)$.
\end{prop}
\dem The equality $\Ker \Th_C= (\ga-1)\U_\infty \cap \C_\infty$ is \ref{ThU}.
It follows that
$$\Ker ( (\C_\infty)_\Ga \overset {\Th_C} \longrightarrow
 \U^\prime_F ) = \Ker ((\C_\infty)_\Ga \overset {\text{nat}} \longrightarrow
(\U_\infty)_\Ga).$$
The latter has been shown in \ref{KNF} to be the $\ZM_p$
torsion of $(\C_\infty)_\Ga$, isomorphic to $(KN_F)^\Ga$.\par\qed
\begin{cor}\label{malbe} Using \ref{ThC} we define the map $\al$
by the commutative triangle~:
\noindent\centerline{\xymatrix {(\C_\infty)_\Ga \ar[r]^-{\Th_C}
\ar@{->>}[d]& \U^\prime_F \\
(\C_\infty)_\Ga/(KN_F)^\Ga \ar[ur]_-\al & , }}
and the map $\be\colon (\C_\infty)_\Ga/
(KN_F)^\Ga \longrightarrow \U^\prime_F/\U_F$ by composing $\al$ with
the natural projection $\U^\prime_F \twoheadrightarrow
\U_F^\prime/\U_F$. They both depend on $\ga$ and are injective.\end{cor}
\dem The map $\al$ is well defined and injective by \ref{ThC}. For $\be$,
the image of $\al$ lies inside $\Uha^\prime_F$, and $\Uha^\prime_F
\cap\U_F=\{1\}$
by lemma \ref{cdcl}. Hence $\be$ is injective.\par \qed
\begin{defi} Put $\pucod_F:=\im \Th_C = \im \al$
and call it the subgroup of {\rm modified circular $S$-units of $F$} (the
terminology will be justified in \ref{Ku=sol} below).
Solomon's ``$S$-unit'' is thus a particular modified circular
$S$-unit.\end{defi}
Obviously $\pucod_F \subset (\U^\prime_\infty)_\Ga \subset \Uha^\prime_F
\subset \U^\prime_F$. By $\be$ in \ref{malbe} we have $\pucod_F =\Im
\al \simeq \Im \be$, and we may also
consider $\pucod_F$ as a subgroup of $\U^\prime_F/\U_F$. The
distinction between $\Im \al$ and $\Im \be$ will always be made
clear by the context.
\begin{prop}\label{indpucod} $\Uha^\prime_F/\pucod_F$ is the
$\ZM_p$-torsion of $\U^\prime_F/\pucod_F$, and its order is
 $$ (\Uha^\prime_F : \pucod_F)= \kappa_F
 \#((X^\prime_\infty)^+)_\Ga\ .$$\end{prop}
\dem Because of the validity of Gross' conjecture, $\Uha^\prime_F$
has $\ZM_p$-rank $(r_1+r_2)$, hence $\Uha^\prime_F/\pucod_F$ is
$\ZM_p$-torsion by \ref{ThC}. Besides, $\U^\prime_F/\Uha^\prime_F$
is $\ZM_p$-torsion free by Sinnott's exact sequence, hence the
first part of the proposition. Let us compute the index
$(\Uha^\prime_F :\pucod_F)= (\Uha^\prime_F :
(\U^\prime_\infty)_\Ga ) ((\U^\prime_\infty)_\Ga : \pucod_F)$.
Using \ref{ThC} we have an isomorphism $(\C_\infty)_\Ga/
(KN_F)^\Ga$ $\overset \sim \longrightarrow \pucod_F$. Using the
snake sequence of the proof of \ref{KNF} $(ii)$ (viz. applying the
snake lemma to multiplication by $(\ga-1)$) we get an isomorphism
$(\U^\prime_\infty)_\Ga / \pucod_F \simeq
(\U_\infty/\C_\infty)_\Ga$. By a classical formula (Herbrand's
quotient in Iwasawa theory), the order of the right hand side is
$p$-adically equivalent to $\kappa_F G(0)$, where $G(T)$ is the
common characteristic series of $(X^\prime_\infty)^+$,
$X^+_\infty$ and $\U_\infty/\C_\infty$. In other terms
$((\U_\infty^\prime)_\Ga : \pucod_F)\overset p \sim \kappa_F \#
((X^\prime_\infty)^+)_\Ga / \# ((X^\prime_\infty)^+)^\Ga$. Note
that by Gross' conjecture, $((X^\prime_\infty)^-)^\Ga$ is finite,
hence null because $ (X^\prime_\infty)^-$ has no non-trivial
finite submodule (see \cite{I73}).
It follows that $((X^\prime_\infty)^+)^\Ga =
(X^\prime_\infty)^\Ga$. As for the quotient
$\Uha^\prime_F/(\U^\prime_\infty)_\Ga$, it is known to be
isomorphic to $(X^\prime_\infty)^\Ga$ (see \cite{Ku72}, 7.5). The
proof of the proposition is complete.\par\qed
\subsection{Regulators}
From now on, we impose the additional condition that $F$ be
{\sl totally real}, in order to get regulator formulae. Note that in this
case, all modules $\pucod_F, \ (\U_\infty^\prime)_\Ga, \ \Uha^\prime_F,$
and $\U^\prime_F/\U_F$ have the same $\ZM_p$-rank $r_1=[F:\QM]$. Note
also that in the totally real case, Leopoldt's conjecture (i.e. the
finiteness of $(\mathfrak X_\infty)^\Ga$, see e.g. \cite{I73} or
\cite{Wi85}) implies what we called
Gross' conjecture (i.e. the finiteness of
$(X_\infty^\prime)^\Ga$) in \S 1.2. However we prefer to keep the terminology
``Gross' conjecture'' because in general the module $\mathfrak
{X}_\infty$ and $X^\pr_\infty$ are not of the same nature : they are
mutually in Spiegelung.
Recall that we have injective homomorphisms (exact sequence ($5\ \la$) and
lemma \ref {cdcl}) between $\ZM_p$-lattices of rank $(r_1-1)$~:
\par
\noindent\centerline{\xymatrix{\U_F \ar@{^{(}->}[r] &
\U^\prime_F /\Uha^\prime_F \ar@{^{(}->}[r]^-\la & p I(S) }}
\par
\noindent We consider $p I(S)$ as a submodule of $I(S)$ and we define
{\sl regulators~:}
\begin{defi}\label{regs} Since the abelian field $F$ satisfies
 both conjectures
of Leopoldt and Gross, we can define
\begin{itemize}
\item[$(i)$] $\rgr_F$
 as the index of $\la(\U_F^\prime/\Uha_F^\prime)$ inside
$I(S)$
\item[$(ii)$]  $\rle_F$ as the index of $\la(\U_F)$ inside
$I(S)$
\end{itemize}\end{defi}

\bigskip

\noindent {\it Remarks}
\begin{itemize}
\item[$(i)$] The exact sequence $(1)$ in \S 1.1 and the isomorphism
$\UC_F \simeq \bigoplus_{v\in S} (1+p\ZM_p)$ (because $p$ is totally
split) show immediately that our index $\rle_F$ is $p$-adically
equivalent to the classical Leopoldt $p$-adic regulator.
Note that we are obliged to work inside $I(S)$ (and not $pI(S)$) because
of Leopoldt's definition.

\item[$(ii)$] The exact sequence $(5)$ in \S 1.2 and the isomorphism
$\oplus_{v\in S} \Fba_v^\times/\Fha_v^\times \simeq \oplus_{v\in
S}(1+p\ZM_p)$ also show that, by taking $p$-adic logarithms, our index
$\rgr_F$ is $p$-adically equivalent to a determinant which is the real
analog of Gross' ``imaginary'' regulator as defined in \cite{FGS}.

\item[$(iii)$] Clearly, $\rgr_F/\rle_F$
is an integer equal to $(\U^\prime_F/\Uha_F^\prime
 : (\U_F\oplus \Uha_F^\prime)/\Uha^\prime_F)= ( \U^\prime_F : \U_F \oplus
\Uha^\prime_F)$.
\end{itemize}

\begin{theo}\label{ind}
Let $F$ be an abelian, totally real number field, and let $p$ be
totally split in $F$.
Let $h^\prime_F$ be the $S$-class number of $F$. Then :
\begin{itemize}
 \item [{\rm $(i)$}] $\#\tor_{\eZ_p} (\U^\prime_F/\widetilde
C_F^{\prime\prime})
=(\widehat U_F^\prime:\widetilde
C_F^{\prime\prime})\overset p \sim\kappa_F \
 h^\prime_F R_F^{\text{Gross}}p^{1-r_1}$
\item [{\rm $(ii)$}] $(\U^\prime_F/\U_F:\widetilde
C_F^{\prime\prime})
=(\U^\prime_F:\widetilde C_F^{\prime\prime}\oplus \U_F)\overset p
\sim\kappa_F \  h_F^\prime
R_F^{\text{Leop}}p^{1-r_1}$
\end{itemize}
(the sign $\overset p \sim$ means $p$-adic equivalence).
\end{theo}
\dem
The first equality in $(i)$ has been shown in \ref{indpucod}. It remains to
compute the $p$-adic valuation of $\#(X^\prime_\infty)_\Ga$. From Sinnott's
exact sequence ($5\ \la$), we get : $\#(X^\prime_\infty)_\Ga \overset
p \sim h^\prime_F \# \coker \la $. This shows $(i)$, and $(ii)$ follows
by the above calculation of $\rgr_F/\rle_F$.\par\qed
\begin{cor}\label{indc} Let $h_F$ be the class number of $F$. Then
$(\eZ_p[S]:\mu(\pucod_F))\overset p \sim \kappa_F
h_F \rle_F p^{1-r_1}$.
\end{cor}
\dem
Recall that $\pucod_F$ is embedded in $\U^\prime_F/\U_F$ by \ref{malbe},
and the map $\mu$ takes place in the exact sequence~:
\par
\noindent \hfill \xymatrix {0\ar[r]& \U^\prime_F/\U_F
\ar[r]^-\mu & \ZM_p [S]
\ar[r]^-{\mathrm {Artin}} &
X_F\ar[r]& X^\prime_F\ar[r]&0}
\hfill ($3\ \mu$).
\par
\noindent Together with \ref{ind} $(ii)$, this shows the corollary.\par\qed
\begin{itemize}
\item[] {\it Remark~:} In \cite{So92}, definition 4.1,
Solomon introduces the Galois
modu\-le $\KC(F)$ generated by all Solomon elements attached to all subfields
of $F$ distinct from $\QM$, and he shows, in the semi-simple case, an index
formula analogous to that of \ref{indc} (but he takes quotients by norms;
see \cite{So92}, proposition 4.3). To compare $\overline{\KC(F)}$ and
$\pucod_F$, see theorem \ref{gen} below.
 Corollary \ref{indc} is clearly a strengthening of Solomon's result.
It gives the most general estimation of the size of the modified circular
$S$-units. Note also that the regulator formula in \ref{indc} bears
a resemblance with Leopoldt's formula giving the residue at $1$ of
the $p$-adic zeta function of $F$. This must (and will) be explained
(see \S \ref{anni}).
\end{itemize}
\subsection{\bf Kuz$'$min's modified circular $S$-units}
Let us give now another description of the modified circular
$S$-units. For simplicity, we stick to the hypotheses of theorem
\ref{ind}. We have seen that for the purpose of descent and
co-descent in Iwasawa theory, it is often more convenient to use
the $S$-units $\U_\infty^\prime$ instead of the units $\U_\infty$
(for example : the natural map $(\U^\pr_\infty)_\Ga
\longrightarrow \U^\pr_F$ is injective by Kuz$'$min's theorem,
here labelled \ref{codes}, whereas $(\U_\infty)_\Ga
\longrightarrow \U_F$ is never injective). It is also natural to
introduce the group $C^\prime_n$ of {\sl circular $S$-units} of
$F_n$ in the sense of Sinnott, which is the intersection of
$U^\prime_n$ with the {\sl circular numbers} of $F_n$ (for details
see \S\ref{cycnum}), and put $\C^\prime_n = C^\prime_n\otimes
\ZM_p$, $\C^\prime_\infty=\limpro \C^\prime_n$. The drawback is
that the $\La$-torsion modu\-le
$\U^\prime_\infty/\C^\prime_\infty$ has no longer finite
co-invariants (contrary to $\U_\infty/\C_\infty$). To get smoother
descent, Kuz$'$min has introduced --without any splitting
hypothesis on $p$-- a certain module of {\sl modified circular
$S$-units at infinite level} which has been studied at length in
\cite{Ku96} and \cite{BN1} (but note that \cite{Ku96} uses
circular units in the sense of Washington -- according to the
terminology of \cite{KN95}-- and \cite{BN1} in the sense of
Sinnott). In our particular situation here ($p$ totally split),
Kuz$'$min's definition can be adapted and much simplified~: let
$\C^{\prime\prime}_\infty$ be the $\La$-submodule of
$\U^\prime_\infty$ such that
$(\U^\prime_\infty/\C^\prime_\infty)^\Ga =
\C^{\prime\prime}_\infty /\C^\prime_\infty$. The obvious guess is
then the good one~:
\begin{prop}\label{Ku=sol} With the hypotheses of theorem \ref{ind},
the inclusion $\C^{\prime\prime}_\infty\subset\U^\prime_\infty$ induces,
by taking co-invariants, an isomorphism
$$\pucod_F\simeq (\C^{\prime\prime}_\infty)_\Ga/(KN_F)^\Ga$$
\end{prop}
\dem The proof will proceed in two steps :
\begin{lem}\label{deprasa} Multiplication by $(\ga-1)$ induces isomorphisms
of $\La$-modules $\C^{\prime\prime}_\infty \overset \sim
\longrightarrow \C_\infty$ and $ \U^\prime_\infty /
\C^{\prime\prime}_\infty\cong \U_\infty /\C_\infty$.  \end{lem}
\dem Let $\tau$ be the inverse of the isomorphism
$\U^\prime_\infty \overset {\overset {\ga-1} \sim} \longrightarrow
(\U^\prime_\infty)^{(0)}$ of theorem \ref{iso}. If $\al_\infty\in
\tau(\C_\infty)$, then $(\ga-1) \al_\infty \in \C_\infty \subset
\C^\prime_\infty$, hence $\al_\infty + \C^\prime_\infty \in
(\U^\prime_\infty/\C^\prime_\infty)^\Ga$ and $\al_\infty \in
\C^{\prime\prime}_\infty$. Conversely, if $\al_\infty\in
\C_\infty^{\prime\prime}$,
 then $(\ga-1)\al_\infty\in \C^\prime_\infty$
by definition. But we have seen in the proof of \ref{dtr} that
$(\ga -1)\al_\infty \in \U_\infty$, hence $(\ga-1) \al_\infty\in \C_\infty$
 i.e. $\al_\infty \in \tau(\C_\infty)$.
(NB : we did not need to suppose
$F$ totally real).\par\qed\par
We now complete the proof of \ref{Ku=sol}. By \ref{deprasa} and the
construction of the map $\Th$, we see that $\pucod_F$ is the image of
$\C^{\prime\prime}_\infty$ by the natural map
$\U^\prime_\infty \rightarrow (\U^\prime_\infty)_\Ga \hookrightarrow
 \U_F^\prime\ .$ It remains only to show that the inclusion
$\C^{\prime\prime}_\infty \subset \U^\prime_\infty$ induces
an injection $(\C^{\prime\prime}_\infty)_\Ga/(KN_F)^\Ga
 \hookrightarrow
 (\U^\prime_\infty)_\Ga$. This follows from \ref{deprasa} and \ref{ThC}.
\par\qed\par
\section {A bridge over troubled water}\label{bridge}
From now on, unless otherwise stated, $F$ will be an {\sl
 abelian number field, totally real and $p$ will be totally split
in $F$}, and let $G=\Gal (F/\QM)$. Our goal is to compare the groups of
circular units and modified circular $S$-units, $\C_F$ and
$\pucod_F$.
\subsection {Statement of the problem}
By construction $\pucod_F \simeq ((\C_\infty)_\Ga)/(KN_F)^\Ga$
(proposition \ref{Ku=sol}) but, as noticed in the introduction,
the natural map $(\C_\infty)_\Ga \longrightarrow \C_F$ gives no
information. The idea is to replace it by the map $(\C_\infty
)_\Ga \longrightarrow \U^\prime_F$ derived from $\Th_C$ of
\ref{dtr}, and to compare its image $\pucod_F$ with the group of
cyclotomic units $\C_F$ inside $\U^\pr_F$. But staying in $\U^\pr_F$
sheds no new light.  Our strategy will be to consider
$\pucod_F$ inside $\U^\pr_F/\U_F$ and $\C_F$ inside
$\U^\pr_F/\Uha^\pr_F$, and the to ``lay a bridge over troubled
water'' between the exact sequences ($3\ \mu$) and ($5\ \la$). As
noticed in \S\ref{SEs}, the water is really troubled, because
there is a priori no natural map between $\U^\prime_F/\U_F$ and
$\U^\prime_F/\widehat U^\prime_F$. For the sake of clarity, let us
give first an abstract of the construction of the bridge~:
\begin{itemize}
\item[--] embed $\C_F$ in $\U^\prime_F/\Uha^\prime_F$ by $\ref{cdcl}$
 and $\pucod_F$ in $\U^\prime_F/\U_F$ by \ref{malbe}
\item[--] then embed $\C_F$ in $I(S)$ via $\la$ and $\pucod_F$ in
$\ZM_p[S]$ via $\mu$
\item[--] use, as in \cite{So92}, Coleman's theory of power series
to compute the semi-local module $\mu (\pucod_F)$
\item[--] this computation suggests to introduce an auxiliary module
$\Dp_F$ of ``circular numbers'', which contains naturally $\C_F$ and is
sent naturally to $\pucod_F$
\item[--] the composite map $\C_F \longrightarrow \Dp_F \longrightarrow
\pucod_F$ gives what we want. Its kernel and cokernel are under control.
\end{itemize}
\subsection{Explicit semi-local description of $\mu\left (\widetilde
C_F^{\prime\prime}\right )$}
We now review --and simplify-- the proof of the main result
of Solomon (see \cite{So92}, theorem 2.1).
Let $L$ be a local $p$-adic field, i.e. a finite extension of $\QM_p$.
As previously, $\Lti$ (resp. $\U(L)$) will denote the $p$-completion
of $L^\times$ (resp of the units), and $\Ltim \infty :=\limpro \Ltim n, \
\U_\infty(L) := \limpro \U(L_n)$ when going up the cyclotomic
$\ZM_p$-extension $L_\infty = \bigcup_{n\geq 0} L_n $ of $L$.
\begin{lem}[local analog of \ref{iso} and \ref{dtr}]\label{isol}
Suppose that $L/\QM_p$ is tamely ra\-mified.
Fix a topological generator $\ga$ of $\Gal(L_\infty/L)$ and let
$N$ be the norm map of $L/\QM_p$. Then multiplication by $(\ga-1)$ induces
an isomorphism
$$\Qtim \infty \overset{ \ga -1} {\overset \sim \longrightarrow}
(\Qtim \infty )^{(0)} = N(\U_\infty(L)).$$
\end{lem}
\dem The first isomorphism (multiplying by $T$) is proved exactly in
the same way as in
 \ref{iso}.
The second equality comes from the tameness
assumption and the surjectivity of the norm map in this case.
\par\qed\par
We consider the special case $L=\eQ_p(\zeta_p)$
and we want to make explicit the inverse of the isomorphism of \ref{isol},
say $\tau_p\colon  N(\U_\infty
(L))\overset \sim \longrightarrow \Qtim \infty $. Let us start from
the exact sequence of $\La$-modules
\xymatrix@=20pt{0\ar[r]& \U_\infty(L)\ar[r] &\Ltim \infty
\ar[r]^-{v_\infty} & \ZM_p \ar[r] & 0}, where
the valuation $v_\infty$ is defined using a choice of
norm coherent uniformizing elements $\pi_\infty=(\pi_n)$.
This exact sequence in general does not split, but never\-theless,
every element $x_\infty\in\Ltim \infty $ may be uniquely written
$x_\infty = z_\infty \pi_\infty^{v_\infty (x_\infty)},$ $z_\infty \in
\U_\infty(L)$. In particular, take
$x_\infty = \tau_p (N(u_\infty))\in \Qtim \infty
\subset \Ltim \infty $. By Coleman's theory over
 $\eQ_p(\zeta_p)$ we can associate uniquely to $u_\infty$
(resp. $z_\infty$) formal power series $g_{u_\infty} (T)$
(resp. $g_{z_\infty} (T)$) in $\Lambda = \eZ_p[[T]]$, $T=\ga -1$.
Let $c$ denotes the topological generator of $1+p \ZM_p$ which corresponds
to $\ga$ by the isomorphism $\ZM_p^\times \cong \Gal(\QM_p(\ze_p)/\QM_p)$.
Then
\begin{lem}[Solomon's lemma]\label{so3.1} (see \cite{So92}, Theorem
3.1).
\par
\noindent Let $\tau_p(N(u_\infty)) = z_\infty \pi_\infty^{a(u_\infty)}
$. We have :
$$\prod_{\omega \in \mu_{p-1}} g_{u_\infty} ((1+T)^\omega -1) =
\left (\frac {(1+T)^c -1} T \right )^{a(u_\infty)} \frac {
g_{z_\infty} ((1+T)^c -1)} {g_{z_\infty} (T) }\ \in \Lambda .$$
In particular $a(u_\infty)=(p-1) \log_c (g_{u_\infty}
(0))$, where $\log_c (\cdot) = \log_p(\cdot)/\log_p(c)$, and $\log_p$
is the Iwasawa logarithm  characterized by $\log_p(p)=0$.
\end{lem}
\dem This is an immediate consequence of the construction of
$\tau_p$. The left hand side product is the norm of Coleman's
power series $g_{u_\infty} (T)$. The right hand side is just obtained
by applying $(\ga-1)$ to $(z_\infty
\pi_\infty^{a(u_\infty)})$.\par\qed\par Let us now come back to
our number field $F$ and give an explicit description of the
semi-local module $\mu(\widetilde C_F^{\prime\prime})\subset
\ZM_p[S] =\oplus_{v\mid p} \eZ_p v $. For any subfield $M$ of $F$,
let us denote by $m$ the conductor of $M$ (recall that $p\nmid
m$). Consider the norm coherent sequence in the $\ZM_p$-extension
of $M$~:
$$\ep_{M,\infty}=\left ( N_{\eQ\left (\ze_{m p^{n+1}}\right )/M_n}
(1-\zeta_m^{\sigma^{-n}}
\zeta_{p^{n+1}})\right )_{n\geq 0}$$
where $\sigma\in \Gal (\eQ(\zeta_{f_M})/\eQ) $ is the Frobenius of
$p$. For $M\neq \QM$, $\ep_{M,\infty}$ is clearly an element of $\C_\infty(F)$.
A system of Galois generators of $\C_\infty(F)$ is given by Greither's lemma
(\cite{Grei92}, lemma 2.3)~:
\begin{lem}\label{genCG}
The elements $\ep_{M,\infty}$
(for $\eQ \subsetneqq M \subset F$),
together with $(\gamma-1) \ep_{\eQ,\infty}$, form a system of
$\Lambda[G]$-generators of $\C_\infty(F)$.
\end{lem}
By applying $\mu\circ\Theta$ to this system of generators we
obtain a system of $\eZ_p[G]$-generators of
$\mu (\widetilde
C_F^{\prime\prime})$, consisting of the elements
$\sum_v a(\imath_v(\ep_{M,\infty})) v$, $\QM \subsetneqq M
\subset F$, $v\in S$ and of $\mu(p)=(1,\ldots,1)$
(here $\imath_v$ is the embedding of $F$ in $F_v$ and $a(\cdot)$ is as
in \ref{so3.1}). We are then led by lemma \ref{so3.1} to evaluate at $0$
the power series $g_{\ep_{M,\infty}}(T)$.
Once a choice of $\pi_\infty$ has been made, these power series are uniquely
 determined. It is possible to choose a sequence
$\pi_\infty$ such that $g_{\ep_{M,\infty}}(T)$ may be easily written down
(actually $g_{\ep_{M,\infty}}(T)$ will turn out to be a polynomial). It then
gives $g_{\ep_{M,\infty}}(0)=\imath_v(\ep_M)$, where for all $M$, of
conductor $m$ say,
$\ep_M$ is the cyclotomic number $\ep_M:=N_{\QM\left (\zeta_m\right )/M}
\left (1-\zeta_m\right )$. Applying lemma \ref{so3.1} one gets
$\mu(\Th_C(\ep_{M,\infty}))=(\log_c(\imath_v(\ep_M)))_{v\in S}$.
(for details see pp. 343
and 344 of \cite{So92}). This proves the :
\begin{theo}\label{gen}
The $\ZM_p[G]$-sub\-mo\-du\-le $\mu(\widetilde
C_F^{\prime\prime})$ of $\eZ_p[S]$ is generated
by all the
$(\log_c(\imath_v(\ep_M)))_{v\in S}$ for $\eQ\subsetneqq M
\subset F$ together with $(1,\ldots,1)$, where
$\log_c(x)$ denotes $\log_p(x)/\log_p(c)$.\par\qed
\end{theo}
\noindent (compare with \cite{So92}, theorem 2.1;
proposition 4.2)
\begin{itemize}
\item[] {\it Remark~:} Because of the injectivity of $\mu$, this shows
the precise relationship between our $\pucod_F$ and Solomon's $\overline
{\KC(F)}$ (\cite{So92}, definition 4.1) :
$$\pucod_F=\langle \KC(F), p
\rangle_{\ZM_p}$$
\end{itemize}
\subsection {Cyclotomic numbers}\label{cycnum}
Let $\Dnb_F$ be the subgroup of cyclotomic numbers of $F^\times$
in the sense of Sinnott. More precisely, $\Dnb_F$ is generated by the elements
 $\ep_M$ (as defined in \ref{gen}) and their Galois conjugates.
These elements are actually local units except possibly at the primes
of $F$ which are ramified above $\QM$. Let us denote by
$\text{Ram}(F/\eQ)$ the set consisting of those primes.
For technical reasons we'll consider $\Dnb_F$ as a subgroup
of the group $U_F(R)$ of $R$-units of $F$, where
$R=\text{Ram}(F/\eQ)\bigcup S\bigcup \{ \text {primes of}\ F $
$\text{dividing}\ (1+p) \}$, and denote by $\D_F$
the closure of $\Dnb_F$ in $\U_F(R)=U_F(R)\otimes \eZ_p$.
We also introduce the group $\Dnbp_F= \langle \Dnb_F, (1+p) \rangle$
and its closure $\Dp_F$.
For any finite set of primes $T$ containing $S$, the ``$T$-analog''
of Sinnott's exact sequence is valid, viz. we have~:
\par
\noindent{\xymatrix{ \U_F(T)
\ar[r]^-{\delta_T}  & \bigoplus_{v\in T}  \overline {F_v^\times}/
\widehat {F^\times_v}\ar[r]^-{\mathrm{Artin}}&
 \Gal (L^\prime_0/F)\ar[r]& X_F(T) \ar[r]  & 0}}  \hfill $(6)$
\par
\noindent Here $\widehat {F^\times_v}$,
the {\sl local cyclotomic norms}, is the
group of norms in the $\eZ_p$-extension $F_{v,\infty}/F_v$
(when $v\notin S$, $\widehat {F^\times_v}$ is the torsion of
$\overline {F_v^\times}$, i.e.
the local $p$-primary roots of unity).
The exactness of $(6)$ is obvious from $(4)$, because
$L_\infty^\prime$, hence $L_0^\prime$, are independent of $T\supset S$.
By Leopoldt's conjecture, $\U_F(T)$ is embedded in $\oplus_{v\in T}
\overline {F^\times_v}$, and $\Ker \de_T= \U_F(T) \cap \oplus_{v\in T}
\widehat {F^\times_v}$. It follows at once from the definition of
$\widehat {F^\times_v}$ that $ker \de_T = \Ker \de_S = \Uha^\prime_F$,
independently of $T\supset S$.
\begin{itemize}
\item[] {\it Remark~:} Jaulent's presentation of Sinnott's exact sequence
(in \cite {J94}) enlarges $T$ to the set of {\it all} primes of $F$.
This has the advantage to dispense with Leopoldt's conjecture.
\end{itemize}
\begin{theo}\label{epimon}
\ \par
\begin{itemize}
\item [{\rm $(i)$}] We have a canonical epimorphism
$\Dp_F \twoheadrightarrow
 \widetilde C_F^{\prime\prime}$.
\item [{\rm $(ii)$}] This epimorphism induces a monomorphism
$\varphi\colon \C_F \hookrightarrow
\ _{\Tr}  (\widetilde C_F^{\prime\prime}) $. Here and from now on,
 $_{\Tr}(\cdot)$ will denote the kernel of the algebraic trace ($=$
action of the
trace element in $\ZM_p[G]$).
\end{itemize}
\end{theo}
\dem
$(i)$ Let $\la_R$ be the composite map $(\oplus_{v\mid p} \log_p)\circ
\text{pr}\circ \de_R$, where pr denotes the projection $\oplus_{v\in R}
\overline {F^\times_v}/\widehat {F^\times_v} \twoheadrightarrow
\oplus_{v\in S} \overline {F^\times_v}/\widehat {F^\times_v}$.
By ($3\ \mu$) and ($5\ \la$), we have a commutative square :
\par
\centerline {\xymatrix {\widetilde C_F^{\prime\prime} \ar[r]^-\mu  &
\eZ_p[S] \\
\Dp_F \ar[r]^-{\la_R}  \ar@{-->}[u]^-\phi & p\eZ_p[S]
\ar[u]^-{\frac 1 p} }}
\par
\noindent where the dotted arrow $\phi$ exists because of \ref{gen} and the
injectivity of $\mu$. The description of the generators in \ref{gen} also
ensures the surjectivity of $\phi$.
\par
\noindent $(ii)$ By definition of $\la_R$, we have another commutative
square
\par
\centerline {\xymatrix{ \Dp_F \ar[r]^-{\lambda_R}  &
p\eZ_p[S] \\ \C_F \ar@{^{(}->}[r]^-{\lambda}  \ar@{^{(}->}[u]  & p\eZ_p[S]
\ar@{=}[u] }}
\begin{itemize} \item[{}] {\it Remark~:}
Actually the image of $\la$ is contained in
$p I(S)$ because the norm of a $(p)$-unit is a power of $p$
(up to a sign). But this holds no longer true
for $T$-units, $T\supsetneqq S $.
\end{itemize}
\noindent Gluing together the two commutative squares, we obtain
still another~:
\par
\centerline{\xymatrix@C=15pt {_{\Tr}\left (\widetilde C_F^{\prime\prime}\right
)  \ar@{^{(}->}[r]^-\mu  & I(S) & \\
\C_F \ar@{^{(}->}[r]^-{\lambda}  \ar[u]^-\varphi  & pI(S)
\ar[u]_-{\frac 1 p}^-\rbag, &\ }}
\noindent
where $\varphi$ is the restriction of $\phi$ to $\C_F$. The injectivity of
$\varphi$ follows from that of the three other sides of the square.
\par\qed\par
The monomorphism $\varphi\colon \C_F \hookrightarrow {}_{\Tr}(\pucod_F)$
 allows us to compare these two $\ZM_p$\-lattices in terms of ramification
and ``structural constants''. Recall {\sl Sinnott's index formula}~:
$(U_F:C_F)=c_F h_F$, where $c_F$ is a rational constant whose explicit
definition does not involve the class group of $F$, and has been
extensively studied in \cite{Si78}, \cite{Si80}, \cite{Ku96}, \cite{BN1},
etc ... Recall that the cyclotomic constant $c$ was introduced in \ref{so3.1},
and $\kappa_F$ in \ref{MFord}.
\begin{theo} Let $p^b\eZ_p$ be the ideal of $\eZ_p$
generated by $[F:\eQ]$ and by the numbers $[F:M_l]\log_c(l)$ where
$M_l:= F\cap \eQ(\zeta_{l^\infty})$ and $l$
runs through all prime divisors of the conductor of $F$
such that $M_l\neq \eQ$.
Then $(_{\Tr}(\widetilde C_F^{\prime\prime}):\varphi(\C_F))$ is
 $p$-adically equivalent to
$ p^b c_F/\kappa_F$.\end{theo}
\dem By definition of $\varphi$, we have a commutative diagram
\par
\centerline {\xymatrix {
0\ar[r] &\C_F \ar[r]^-\varphi \ar@{^{(}->}[d] & \ _{\Tr}(\widetilde
C^{\prime\prime}_F) \ar[r] \ar@{^{(}->}[d]_-{\mu\log_p(c)}  &\coker
\varphi \ar[r]\ar[d] &0  \\
0\ar[r]  & \U_F \ar[r]^-\lambda  &p I(S)\ar[r] &\coker
\lambda \ar[r] & 0 }}
\par
\noindent and since $\log_p(c)\sim p$, we deduce  equalities
$$\left ( I(S):\mu(\
_{\Tr}(\widetilde C_F^{\prime\prime}))\right )=\left (pI(S):\im
\log_p(c)\mu \right )=\left ( \U_F:\C_F\right )
 \frac {\# \coker \lambda} {\#\coker \varphi}\ .$$
But $\#\coker\lambda \overset p \sim
R_F^{\text{Leop}}p^{1-[F:\eQ]}$ by \ref{regs}, and
$(\U_F:\C_F)\overset p \sim c_F h_F$ by Sinnott's formula.
It remains to compute the index
$(I(S):\mu(\ _{\Tr}(\widetilde C^{\prime\prime}_F)))$. Using
the commutative diagram
\par
\centerline {  \xymatrix {
    0\ar[r] &  \ _{\Tr}(\pucod_F)\ar@{^{(}->}[d]_-\mu \ar[r] &
 \pucod_F \ar[r]^-{\Tr}  \ar@{^{(}->}[d]_-\mu  &
{\Tr}(\pucod_F) \ar[r] \ar[d]_-\psi  & 0 \\
0\ar[r] & I(S) \ar[r] &\eZ_p[S]\ar[r]^-{\Tr} & \eZ_p
\ar[r] & 0               }}
\noindent it appears that we just have to study the vertical map
$\psi$. A priori we know that ${\Tr}(\pucod_F)$ has $\eZ_p$-rank one.
But here we may identify the algebraic norm ${\Tr}$ with the arithmetic
norm $N$ from $F$ down to $\eQ$ (recall
 that we are dealing with $S$-units),
so that ${\Tr}(\pucod_F)$ may be identified with a submodule of
$\U^\prime_\eQ=p^{\eZ_p}$; in particular it has no
$\eZ_p$-torsion. This implies that
$\Ker \psi = 0$. Concerning
the image of $\psi$, the commutativity of the diagram
shows that $\im \psi$ is the ideal of
$\eZ_p$ generated by the elements
$\sum_{v\in S} \log_c(\imath_v(\ep_M))=\log_c(N(\ep_M))$,
together with $[F:\eQ]={\Tr}(\mu(p))$.
Even more explicitly, we know that if the conductor $m$ of $M$
is a power of a single prime, $\ell$ say, then
 $N(\ep_M)=\ell^{[F:M]}$, and if $m$ is composite, $N(\ep_M)=1$.
Denoting (as indicated in the statement of the theorem)
$\im \psi = p^b \eZ_p$, we have
$(I(S):\mu(\ _{\Tr}(\pucod_F)))=p^{-b}
(\eZ_p[S]:\mu(\pucod_F))$ and this last quantity is $p$-adically equivalent
to $\kappa_F p^{-b+1-r_1} h_F R_F^{\text{Leop}}$ by
corollary \ref{indc}. The theorem is proved.\par\qed
\begin{itemize}
\item[] {\it Remark~:} In \cite{So94}, Solomon uses a description by
generators
and relations to compare (in our notations) $\C_F$ and $\pucod_F$. In doing
so, he must enlarge these modules to get (in his notations) a surjective
map $\overline\th\colon \overline {\DC_\SC}\twoheadrightarrow \overline
\KC_\SC $ (\cite{So94}, theorem 1) which plays a role analogous
to our map $\phi$ in the proof of \ref{epimon}. The restriction of
$\overline {\th}$ to $\C_F$ coincides with our map $\varphi$, and
theorem 4.2 of \cite{So94} determines its image. This construction
by generators and relations allows one to compute $\ZM_p$-ranks,
but probably not indices nor annihilators.
\end{itemize}
\section{Applications}\label{anni}

\subsection{Annihilation of class groups}\label{solcon}
As announced in the introduction, we are now looking for ideals of
$\ZM_p[G]$ which annihilate the $p$-groups $X_F$ and $\U_F/\C_F$.
Even more precisely, we intend to compute the Fitting ideals of some related
modules. In order to state Solomon's conjecture, let us temporary
relax the hypothesis that $p$ is totally split in the totally real abelian
field $F$.
So now we only assume that
$p$ is a prime number which doesn't ramify in $F$, which in turn may
be a real or complex abelian field.
We fix an embedding
$\imath \colon \eQ^{\text{sep}}
\hookrightarrow  \eQ_p^{\text{sep}}$, and we shall denote
by $\OC$ the topological closure in $\CM_p$ of the image $\imath(O_F)$
of $F$ (a priori $\OC$ is inside $\CM_p$, but it actually lies
in $\QM_p^{\text{sep}}$). For any subfield $M\neq \QM$ of $F$, we define
{\sl Solomon's element} $\sol_M$
as in \cite{So92}, \S 4 :
$$\text{sol}_M:=\frac 1 p \sum_{g\in \Gal (M/\eQ)} \left (
\log_p(\imath(\ep_M^g))\right )g^{-1} \in \OC[\Gal(M/\eQ)].$$
By convention $\sol_\QM =1$. These elements should be considered as real
analogues of the classical Gau\ss\ sums, and the real analog of
Stickelberger's theorem is {\sl Solomon's conjecture} :
\begin{con}[see \cite{So92}, 4.1]\label{conj}
$\text {\rm sol}_F$ annihilates $Cl_F\otimes \OC $.\end{con}
Solomon's conjecture holds true in the semi-simple case (i.e. $p\nmid
[F:\eQ]$), but this does not give anything new, because in this case
it follows from Mazur-\-Wiles' theorem (\cite{So92}, 4.1; see also
\ref{compare} below).
In the general case, let us introduce a modified Solomon's element.
\begin{defi}\label{defSe}  Let $F$ be a totally real abelian
field. Recall that $w\mid p$ is the place corresponding
to the fixed embedding $\imath\colon \QM^{\text{sep}}
\hookrightarrow \QM_p^{\text{sep}}$. We still write
$w$ for the restriction of $w$ to any subfield $M$ of $F$.
\begin{itemize}
\item[{\rm $(i)$}] Let $M$ be any subfield of $F$.
Let $\widetilde{\sol}^F_M$ be any element in $O_{F_w}[\Gal(F/\QM)]$ which
restricts to $\sol_M$ in $O_{M_w}[\Gal(M/\QM)]\subset O_{F_w}[\Gal(M/\QM)]$.
In the commutative ring $O_{F_w}[\Gal(F/\QM)]$, we define the product
$$\sol^F_M :=  \widetilde{\sol}^F_M \ {\Tr}_{F/M}.$$
(with obvious notation for the trace). This product does not
depend on the choice of the lift $\widetilde{\sol}^F_M$.
\item[{\rm $(ii)$}] For any real abelian field $M$ with prime power conductor,
 let $g_M$ be a genera\-tor of the (necessarily cyclic, $M$ being real)
group $\Gal(M/\QM)$. For any subfield $M$ of $F$
we define
$$
\sol_{M,2}^F
:= \left \{ \begin{aligned} &(1-g_M) \sol^F_M \text {if the conductor
of } M \text { is a prime power,}
\\ & \sol^F_M \text{otherwise.}\end{aligned} \right .
$$

\end{itemize}
\end{defi}

We intend now to prove a slightly modified version of
Solomon's conjecture in the
non semi-simple case, with the additional hypothesis that
$p$ is totally split in $F$ (then $\OC =\ZM_p$). It is here that our
functorial approach will pay off, in that it will allow us to apply
techniques ``\`a la Rubin'' to annihilate $X_F$ :
\begin{lem}\label{Rub} Let $N$ be a power of $p$.
Let $V\subset \U_F/(\U_F)^N $ be a Galois
submodule and let $\rh\colon V\longrightarrow
\eZ_p[G]/N \eZ_p[G]$ be any equivariant homomorphism.
Let $\CC\subset U_F$ be the subgroup of "special units" as defined
in {\rm \cite{Ru87}, p. 512}, and write $(\overline \CC)_N$ for  its
image $(\CC (\U_F)^N)/ (\U_F)^N \subset  \U_F/(\U_F)^N$.
Then $\rh(V\cap (\overline \CC)_N)$
annihilates $X_F/N X_F$.
\end{lem}
\dem This lemma is a direct consequence in our
special case of Rubin's theorem 1.3 (see \cite{Ru87}).
The field denoted $K$ in loc. cit. is here equal to $\QM$.
We take for $A$ of loc. cit. the full group $X_F/N X_F$.
Then Rubin's theorem 1.3 shows that $\rh(V\cap \overline \CC)$
annihilates some submodule $A^\prime\subset A$. Using
lemma 1.6 $(ii)$ and the definition of $H_1$ in loc. cit.
it is easily checked that, in our special case
(i.e. $p\neq 2$, the only roots of unity in $F$ are $\pm 1$, and
 every place above $p$ is totally ramified in $F(\ze_N)/F$),
we actually have $A^\prime=A=X_F/N X_F$.\par\qed
\begin{theo}\label{theo}
\par
\noindent  For any totally real abelian field $F$, and any $p$ totally
split in $F$, $\sol_{F,2}^F$ annihilates $X_F$.
\end{theo}
N.B. : If the conductor $f$ of $F$ is divisible by at least two
distinct primes this theorem is Solomon's conjecture \ref{conj},
because by definition the elements $\sol_{F,2}^F$ and $\sol_F$ are
equal. If the conductor is a power of a single prime $\ell$ this
theorem is slightly weaker than Solomon's conjecture. But it is, for
ideal classes, the perfect analogue of an anihilation result
for unit classes that we shall prove later (see theorem \ref{fit}).\\
\smallskip
\dem Let $\rh_1 \colon \U_F \longrightarrow \ZM_p[G]$ be the composite
map $\rh_1 = \eta\circ \frac 1 p \circ \la$, where $\eta$ is the isomorphism
$\eta\colon\eZ_p[S]=\eZ_p[G]\ w
\overset \sim \longrightarrow \eZ_p[G]$, and $w$ is the prime
corresponding to $\imath$. Fix $N$ such that $(X_F/N X_F)=X_F$.
We apply lemma \ref{Rub} by taking $V=\U_F/(\U_F)^N$
and $\rh\colon V\longrightarrow \ZM_p/N\ZM_p[G]$
to be the map induced by $\rh_1$.
Exactly as we defined $\sol_{F,2}^F$ we put :
$$ \ep_{F,2}
:= \left \{ \begin{aligned} &(1-g_F) \ep_F\ \text {if the conductor
of } F \text { is a prime power,}
\\ & \ep_F\ \text{otherwise.}\end{aligned} \right . $$
Then it is a classical fact that $\ep_{F,2}$ is a unit and
the element
$\sol_{F,2}^F$ is nothing else but the image $\rh_1(\ep_{F,2})$.
Further $\ep_{F,2}$ is a special unit :
in the special case that $F$ is the maximal real subfield
of a cyclotomic field this is theorem 2.1 of loc. cit., and
for any abelian field the result for $\ep_{F,2}$ follows from exactly the same
computation. By lemma \ref{Rub} it follows that the class in
$\ZM_p[G]/N$ of $\sol_{F,2}^2$ anihilates $X_F/N
X_F=X_F$. \par\qed\par
\subsection{Fitting ideals of quotients of units}\label{varunits}
Theorem \ref{theo} is still unsatisfactory, in that it uses only
the map $\la$, not the map $\mu$ (notation of \S \ref{SEs}). Let us
cross the bridge built in \ref{epimon}. We'll rather work with ideals than
with elements :
\begin{defi}\label{defSi}
 We define ideals
$ \Sol_1(F)\supset \Sol_2(F) $ of $O_{F_w}[G]$ by giving sets of
gene\-rators~:
$$ \Sol_1(F):=\left \langle  \sol_M^F \vert \QM\subseteq M \subseteq F  \right
\rangle \supset \Sol_2(F):=\left
\langle  \sol_{M,2}^F \vert \QM\subseteq M \subseteq F  \right
\rangle .
$$
\end{defi}

\begin{prop}\label{idS} Suppose that $p$ is totally split in $F$.
Recall that $\eta$ is the isomorphism $\eta\colon\eZ_p[S]=\eZ_p[G]\ w
\overset \sim \longrightarrow \eZ_p[G]$. Then~:
$$ \Sol_1(F)=\eta\circ\frac 1 p \circ\la_R (
 \overline { \text{\rm Cyc}^\prime_F } )\
\text{and} \
\Sol_2(F)=\eta\circ\frac 1 p \circ\la_R(\C_F\oplus (1+p)^{\eZ_p}).$$
\end{prop}
N.B. : The set of places $R$ has been chosen in order that
$\U_F(R)$ contains simultaneously  $ \text{\rm Cyc}^\prime_F $, $\C_F$ and
$(1+p)^{\eZ_p}$, the last two being direct summands.

\dem By definition of the map $\la_R$, we have
$\sol_M^F =\eta\circ\frac 1 p \circ\la_R (\ep_M)$.
Since these elements, together with $(1+p)$ (which is
sent to ${\Tr}_{K/\QM}$, up to a $p$-adic unit)
generate $\overline {\text{\rm Cyc}^\prime_F}$, the first equality is proved.
The second is proved in exactly the same way.\par\qed\par
We are now in a position to prove global annihilation results
for various quotients of units by circular units~:
\begin {theo}\label{fit}
Let $F$ be a totally real abelian field, let $p$ be totally split in $F$.
Then :
\begin{itemize}
\item [{\rm $(i)$}] Let $\UC_n$ be the semi-local module $\UC_n:=\oplus_{v\in
S} U^1_v(F_n)$ and $\UC_\infty := \limpro \UC_n$.
 The embedding $\C_\infty \hookrightarrow \UC_\infty$
induces, by taking co-invariants, an embedding  $\pucod_F\hookrightarrow
(\UC_\infty)_\Gamma$, and $\Sol_1(F)$ is the initial Fitting ideal of
the $\eZ_p[G]$-quotient modu\-le.
\item [{\rm $(ii)$}] Let $\varphi\colon \C_F \hookrightarrow\pucod_F$
be the map defined in \ref{epimon}, and extend it to a map
 $\widetilde\varphi\colon\C_F\oplus (1+p)^{\eZ_p}
\longrightarrow\pucod_F $ by putting $\widetilde \varphi(1+p)=p$.
Then $\mu\circ \widetilde \varphi$ gives an embedding
$\C_F\oplus (1+p)^{\eZ_p} \hookrightarrow (\Cal U_\infty)_\Gamma$, and
$\Sol_2(F)$ is the initial Fitting ideal of the $\eZ_p[G]$-quotient module.
\item [{\rm $(iii)$}] $\Sol_2(F)$ annihilates $\U_F/\C_F$.
\end{itemize}
\end{theo}
\dem
$(i)$ The analog of the exact sequence $(1)$ at infinite level gives an
injection $\UC_\infty/\U_\infty \hookrightarrow \XG_\infty := \limpro
\XG_{F_n}$. Leopoldt's conjecture is known to be equivalent to $\XG_\infty^\Ga
= 0$, hence implies $(\UC_\infty/\U_\infty)^\Ga =0$.
It follows that $(\U_\infty) _\Ga \hookrightarrow
(\UC_\infty)_\Ga$ and that $(\C_\infty)_\Ga \longrightarrow (\UC_\infty)_\Ga$
has the same kernel as $(\C_\infty)_\Ga \longrightarrow (\U_\infty)_\Ga$,
hence the embedding $\pucod_F
\hookrightarrow (\UC_\infty)_\Gamma$.
Let us introduce another semi-local module, namely $\FC_\infty(F)
:= \limpro \oplus_{v\in S} \overline {F_{n,v}^\times}$. We have
seen that multiplication by $(\ga-1)$ induces an isomorphism
$\C^{\prime\prime}_\infty \overset \sim \longrightarrow \C_\infty$
(by \ref{deprasa}) and analogously, $\FC_\infty (F)
\overset \sim \longrightarrow N(\UC_\infty(F(\ze_p)))$ (by \ref{isol}).
But $N(\UC_\infty(F(\ze_p)))=\UC_\infty(F)$ by tame ramification.
In other words, we have a commutative square
\par
\noindent\centerline{\xymatrix{ \C_\infty^{\prime\prime} \ar[r]^-\sim_-{\ga-1}
\ar@{^{(}->}[d] & \C_\infty\ar@{^{(}->}[d] \\
\FC_\infty(F) \ar[r]^-\sim_-{\ga-1} &\UC_\infty(F)=\UC_\infty }}
\par
\noindent By codescent we then get a commutative diagram :
\par
\noindent\centerline{\xymatrix@=18pt{
0\ar[r] & \pucod_F = (\C_\infty^{\prime\prime})_\Ga/(KN_F)^\Ga
\ar[r]^-\mu\ar[d]^-\rbag
 & \FC_\infty(F)_\Ga \cong \ZM_p[S]\ar[d]^-\rbag \ar[r] & \coker \mu
\ar[d]^-\rbag \ar[r] & 0 \\
0 \ar[r]& (\C_\infty)_\Ga/(KN_F)^\Ga
 \ar[r] & (\UC_\infty)_\Ga
\ar[r] & (\UC_\infty/\C_\infty)_\Ga \ar[r] & 0 }}
\par
\noindent which shows that $(\UC_\infty/\C_\infty)_\Ga$ has the same
$\ZM_p[G]$-Fitting ideal as $\coker \mu$. By the isomorphism
$\eta \colon \ZM_p[S] \overset \sim \longrightarrow \ZM_p[G]$,
this last Fitting is nothing but the ideal $\eta \circ \mu (\pucod_F)
= \eta\circ \frac 1 p \circ \la_R (\Dp_F) = \Sol_1(F)$.
\par
\noindent $(ii)$ The commutative diagram
\par
\noindent\centerline{\xymatrix{ \C_F \oplus (1+p)^{\ZM_p} \ar@{^{(}->}
[r]^-{\tilde \varphi}
\ar@{^{(}->}[d] & \pucod_F \ar@{^{(}->}[d]^-\mu \\
\U_F \oplus (1+p)^{\ZM_p} \ar@{^{(}->}[r]^-{\frac 1 p \circ \la}
 &\ZM_p[S] }}
(see the construction of $\varphi$ in \ref{epimon}) and the isomorphism
$\ZM_p[S] \simeq (\UC_\infty)_\Ga$ in $(i)$ above show that
$\mu \circ \tilde \varphi$ embeds $\C_F \oplus (1+p)^{\ZM_p}$ into
$(\UC_\infty)_\Ga$. We conclude by applying $\eta$ as in $(i)$.
\par
\noindent $(iii)$ It is well known that the Fitting ideal of a module
is contained in its annihilator. Here $\Sol_2 (F)$ will annihilate
$(\UC_\infty)_\Ga / \mu\circ\tilde\varphi (\C_F \oplus (1+p)^{\ZM_p} )$,
which contains, by the commutative square in $(ii)$, an isomorphic image
of $\U_F/\C_F$.\par\qed
\begin{cor} (Real analog of Stickelberger's index)
$$(\eZ_p[G]:\Sol_1(F))=(\eZ_p[S]:\mu(\pucod_F))=\kappa_F\#\tor_{\eZ_p}
\mathfrak {X}_F \overset p \sim
\kappa_F h_F\rle_F p^{1-r_1}$$\end{cor}
\dem It has been shown in \ref{fit} $(i)$ that
$(\eZ_p[G]:\Sol_1(F))=(\eZ_p[S]:\mu(\pucod_F))=\#(\UC_\infty/\C_\infty)_\Ga$.
Since $\XG_\infty$ has no non-trivial $\Ga$-invariants
(Leopoldt's conjecture) and has the same characteristic
series as $\UC_\infty/\C_\infty$, we deduce the equivalence :
$\#(\UC_\infty/\C_\infty)_\Ga\overset p \sim \#(\UC_\infty/\C_\infty)^\Ga
\#(\XG_\infty)_\Ga$.
As $F$ is totally real, $(\XG_\infty)_\Ga$ is obviously isomorphic
to $\tor_{\ZM_p} \XG_F$. It remains to compute the order of
$\tor_{\ZM_p} \XG_F$~: this a classical calculation using
$p$-adic $L$-functions (see e.g. \cite{BN1}, 2.6), which is actual\-ly
equivalent to Leopoldt's $p$-adic formula (the
discriminant does not appear here because
$p$ is unramified in $F$).\par\qed
\begin{itemize}
\item[{}] {\it Remark~:} This gives another proof of formula \ref{indc},
at the same time explaining the parenthood between this index formula and
Leopoldt's residue formula.
\end{itemize}
\subsection{Still another exact sequence}\label{zerotriv}
Let us extract from the proof of \ref{fit} the following analog of the exact
sequences $(1)$ and $(2)$ of class field theory~:
\begin{cor}\label{se7}
We have an isomorphism $(\UC_\infty/\C_\infty)_\Ga \simeq
(\oplus_{v\in S} \Fha^\times_v)/\pucod_F$ (notation of \S \ref{Si})
and an exact sequence of finite modules~:
\par
\noindent\centerline{\xymatrix@=20pt{0\ar[r] &
\Uha^\prime_F/\pucod_F \ar[r] & (\oplus_{v\in S} \Fha^\times_v )/\pucod_F
\ar[r] & \tor_{\ZM_p} \XG_F \ar[r] & (X^\prime_\infty)_\Ga \ar[r] & 0 }
\hfill $(7)$}
\end{cor}
\dem It has been shown in the proof of \ref{fit} that
$(\UC_\infty/\C_\infty)_\Ga \simeq (\FC_\infty(F))_\Ga/\pucod_F$, where
$\FC_\infty(F):=\limpro \oplus_{v\in S} \overline {F^\times_{n,v}}$
(notation of \S \ref{Si}). But $(\overline {F^\times_{v,\infty}})_\Ga=
(\overline {F^\times_{v,\infty}})_{\Ga_v}\simeq \Fha^\times_v$ by class field
theory. Thus the isomorphism is proved.
\par
From the exact sequence :
\par
\noindent\centerline{\xymatrix@=19pt {0\ar[r]& \U^\prime_\infty \ar[r] &
\FC_\infty(F) \ar[rr] \ar@{->>}[rd]
& & \XG_\infty \ar[r] & X^\prime_\infty \ar[r] &  0\\
 & & & \DC_\infty \ar@{^{(}->}[ru]& & & }}
(where $\DC_\infty$ denotes the relevant subgroup generated by decomposition
subgroups), we get by co-descent two exact sequences~:
\par
\noindent\centerline{\xymatrix {0\ar[r] & (\U^\prime_\infty)_\Ga \ar[r] &
\oplus_{v\in S} \Fha^\times_v \ar[r] & (\DC_\infty)_\Ga \ar[r] & 0 \quad} and }
\par
\noindent\centerline{\xymatrix {0\ar[r] & (X^\prime_\infty)^\Ga \ar[r] &
(\DC_\infty)_\Ga \ar[r] & (\XG_\infty)_\Ga \simeq \tor_{\ZM_p} \XG_F
 \ar[r] & (X^\prime_\infty)_\Ga \ar[r] & 0 }}
 Putting them together and using Kuz$'$min's isomorphism $\Uha^\prime_F /
(\U_\infty^\prime)_\Ga \simeq (X^\prime_\infty)^\Ga$ (see \cite{Ku72}, 7.5), we
get the exact sequence (which is general)~:
\par
\noindent\xymatrix@=18pt{0\ar[r] & \Uha^\prime_F \ar[r] &
\oplus_{v\in S} \Fha^\times_v \ar[r] & \tor_{\ZM_p} \XG_F
 \ar[r] & (X^\prime_\infty)_\Ga \ar[r] & 0 }, hence also $(7)$. \par\qed
\begin{itemize}
\item[{}] {\it Remark~:}
In the study of the Main Conjecture via circular units inside semi-local
units and Coleman's theory, and even in the semi-simple case, a problem arises
for characters $\psi$ of $G=\Gal(F/\QM)$ such that $\psi(p)=1$, because the
natural co-descent map $(\UC_\infty/\C_\infty)_{\Ga,\psi} \longrightarrow
(\UC_F/\C_F)_\psi$ then gives no information (see e.g. \cite{Gi792}). This
difficulty is related to the phenomenon of trivial zeroes of $p$-adic
$L$-functions. In this context, the more ``sophisticated'' exact sequence
$(7)$ could be viewed as a device to ``bypass trivial zeroes''. We hope
to come back to it in greater detail.
\end{itemize}
\subsection{A criterion for Greenberg's conjecture}\label{greeconj}

Greenberg's conjecture asserts that the orders of the groups $X_n=X_{F_n}$ are
bounded in the cyclotomic $\ZM_p$-extension of a totally real number
field $F$. This conjecture has been
numerically verified by various authors in a huge amount of
special cases. The case when $p$ splits in the base field $F$ is
notably more difficult than the others, probably because of the
phenomenon of trivial zeroes. Our approach here gives a new proof of a
known criterion (see \cite{Ta99}, theorem 1.3; but note that Taya only
needs Leopoldt's conjecture)~:
\begin{theo}\label{greespli} Let $F$ be a real abelian field and $p$ be an odd
prime totally split in $F$. Let $D_n\subset X_n$ be the submodule
generated by (images of) the places of $F_n$ above $p$. Then Greenberg's
conjecture holds true for $F$ and $p$ if and only if $\#
D_n=\#\rm{tor}_{\ZM_p}(\XG_F)$ for every $n$ sufficiently large.
\end{theo}
\dem By definition $D_n$ is the kernel $D_n=\Ker
(X_n\longrightarrow X^\pr_n)$. Hence it may be used to
split the sequence $(3\ \mu)$ in two shorter sequences. Taking
projective limits we get :

\noindent\hfill \xymatrix {0 \ar[r] & \U_\infty \ar[r] &
\U^\pr_\infty \ar[r]^-\mu & \ZM_p[S] \ar[r] & D_\infty \ar[r] &
0}\hfill $(3 \mu \, a)$

\noindent\hfill\xymatrix {0\ar[r]& D_\infty \ar[r] & X_\infty
\ar[r] & X_\infty^\pr \ar[r] & 0} \hfill $(3 \mu\, b)$

\noindent Since $F_\infty/F$ is totally ramified at all places above
$p$, $D_\infty$ is a submodule of $X_\infty^\Ga$ and  therefore is
finite. Again by total ramification, the maps $D_{n+1}
\longrightarrow D_n$ are surjective. So the theorem asserts that
Greenberg's conjecture is equivalent to the equality $\# D_\infty =
\#\rm{tor}_{\ZM_p} (\XG_F)$. From corollary \ref{dtr} we have
$\U^\pr_\infty/\U_\infty=(\U^\pr_\infty)_\Ga$, hence by
$(3 \mu \, a)$ we get the isomorphism
$D_\infty\cong\ZM_p[S]/\mu((\U^\pr_\infty)_\Ga)$. Using
corollary \ref{indc} we deduce the equality :

\noindent \xymatrix {\#D_\infty (\mu((\U^\pr_\infty)_\Ga):\mu
(\pucod_F)) =\#D_\infty
((\U^\pr_\infty)_\Ga:\pucod_F)=\kappa_F\#\rm{tor}_{\ZM_p}(\XG_F)}
\hfill $(\dag)$
\par
Now Greenberg's conjecture is equivalent to the finiteness of
$\U_\infty/\C_\infty$, which in turn is equivalent to the equality
$\U_\infty/\C_\infty=KN_F$. Crossing the bridge, actually applying
lemma \ref{deprasa}, this equality becomes equivalent to the
canonical isomorphism $KN_F\cong
\U^\pr_\infty/\C^{\pr\pr}_\infty$. By Nakayama's lemma, this
isomorphism is equivalent to the equality
$((\U^\pr_\infty)_\Ga:\pucod_F)=\#(KN_F)_\Ga$. Since $KN_F$ is
finite we already have the equalities
$\#(KN_F)_\Ga=\#(KN_F)^\Ga=\kappa_F$. Comparison with the equality
$(\dag)$ shows the theorem.
\par\qed\par
\noindent For numerical applications of criterion \ref{greespli}, see
\cite{Ta99} and the references therein.
\section{Comparison of Fitting ideals}\label{MC}
At this point, a natural question arises : does theorem \ref{fit} give
anything new with respect to the Main Conjecture ? In this section, the
general hypothesis will be~: $F$ is an {\it abelian totally real field}
(and $p$ is not necessarily totally split in $F$),
$G=\Gal(F/\QM)$. To avoid petty technical trouble, let us also
suppose that $F$ {\it is linearly disjoint from }$\QM_\infty$
(which is the case if $p$ is totally split in $F$), in order that
$G$ acts on all the natural modules attached to the fields $F_n$ in the
cyclotomic $\ZM_p$-extension $F_\infty/F$.
\begin{defi} Let $MW(F)$ (for Mazur-Wiles) be the initial Fitting
ideal of $\tor_{\ZM_p} \mathfrak{X}_F$ over $\ZM_p[G]$.
\end{defi}
This ideal annihilates $\tor_{\ZM_p} \mathfrak{X}_F$, hence also
(since $F$ is supposed linearly disjoint from $\QM_\infty$) the $p$-class
group $X_F$, but we do not know a priori if it annihilates $\U_F/\C_F$.
Conversely, if $p$ is totally split in $F$, recall that the ideal $\Sol_2(F)$
annihilates both $\U_F/\C_F$ and $X_F$.
We would like to compare $\Sol_1(F)$ and $MW(F)$.
\par
At infinite level, class-field theory gives an exact sequence of
$\La[G]$-modules~:
\par
\noindent \centerline{\xymatrix { 0\ar[r]& \U_\infty/\C_\infty
\ar[r]& \UC_\infty/\C_\infty \ar[r]& \mathfrak{X}_\infty \ar[r]& X_\infty
\ar[r]& 0}}
\noindent We get by co-descent and functoriality~:
$$MW(F)=\mathrm{Fitt}_{\ZM_p[G]} (\tor_{\ZM_p} \mathfrak{X}_F)
= \text{image in }\ZM_p[G]\text{ of }
\mathrm{Fitt}_{\La[G]} (\mathfrak{X}_\infty)$$
and $$\Sol_1(F)=\mathrm{Fitt}_{\ZM_p[G]} ((\UC_\infty/\C_\infty)_\Ga)
=\text{image in } \ZM_p[G] \text{ of }
\mathrm{Fitt}_{\La[G]}(\UC_\infty/\C_\infty)$$
The Main Conjecture implies that $\UC_\infty/\C_\infty$ and
$\mathfrak{X}_\infty$ have the same characteristic series in $\La$,
but to compare their $\La[G]$-Fitting ideals, more is needed, some
kind of ``equivariant Iwasawa theory'' which does not exist yet
(but there is work in progress, see e.g. \cite{BG00}, \cite{RW02} ...).
Besides, the Fitting ideal of a module is a rather weak invariant, except
if this module is of projective dimension at most $1$ (see the
comments in \cite{Grei98}). For all these reasons, we'll be content
to work in the following setting (which was already that of \cite{BB98},
\cite{BN1}, \cite{Grei98})~:
\begin{nots} Write $G=P\times \De$, where $P$ (resp. $\De$) is the
$p$-part (resp. non-$p$-part) of $G$.
\end{nots}
For any $\QM_p$-irreducible character $\psi$ of $\De$, let
$e_\psi$ be the usual idempotent, and for any $\ZM_p[G]$ or
$\La[G]$-module $M$, let $M_\psi$ be $e_\psi M$. We intend
to compare $MW(F)_\psi$ and $\Sol_1(F)_\psi$ for some non trivial $\psi$.
Let us denote by $\ZM_p(\psi)$ (resp. $\La(\psi)$) the rings
$e_\psi \ZM_p[\De]\simeq \ZM_p[\ch(\De)]$ (resp.
$e_\psi \La[\De]\simeq \La[\ch(\De)]$), where $\ch$ is any
$\QM_p^{\text{sep}}$-irreducible character of $\De$ dividing $\psi$.
\par
We must first show results on projective dimensions, and for that
introduce some more appropriate objects and hypotheses~:
\begin{defi}\label{deflnt} \ \par
\begin{itemize}
\item[$(i)$] Let $T=S\bigcup \mathrm{Ram}(F/\QM)$ (recall
that $S$ is the set of primes above $p$). Slightly abusing language,
we'll keep the notations $S$ and $T$ when going up the cyclotomic
tower. We define
$\mathfrak{X}^T_\infty$ as the
Galois group over $F_\infty$ of the maximal abelian pro-$p$-extension
of $F_\infty$
which is unramified outside $T$.
\item[$(ii)$] Let $\Si$ be the set of places $v\in T-S$ which split
totally in $F(\ze_p)/F$. A $\QM_p$-irreducible character $\psi$ of
$\De$ will be called ``locally non Teichm\"uller'' if either $\Si$ is
empty or the restriction of $\psi$ to the decomposition subgroup
$\De_v$, for any $v\in\Si$, differs from the restriction of the
Teichm\"uller character.
 (Examples of fields for which any
character of $\De$ is locally non Teichm\"uller may be found in
\cite{BB98} and \cite{Grei98}).
\end{itemize}
\end{defi}
\begin{lem}\label{dpX} \ \par
\begin{itemize}
\item [$(i)$]
For any $\QM_p$-irreducible character $\psi$ of $\De$, $\psi\neq 1$, we have
$\text{\rm pd}(\mathfrak {X}^T_\infty)_\psi \leq 1$ over the algebra
$\Lambda(\psi)[P]$.
\item[$(ii)$] For any locally non Teichm\"uller character $\psi\neq 1$
of $\De$, we have $\text{\rm pd}(\mathfrak {X}_\infty)_\psi \leq 1$ over
the algebra $\Lambda(\psi)[P]$.
\end{itemize}
\end{lem}
($(ii)$ generalizes corresponding results of \cite{BB98} and
\cite{Grei98}).
\par

\noindent\dem By \cite{Ng84}, proposition 1.7, there exists a
canonical $\La[G]$-module $\YC_\infty^T$ of projective dimension at most
$1$, such that $\mathfrak {X}_\infty^T$ and $\YC_\infty^T$ take place in
an exact sequence :

\noindent\centerline{\xymatrix{ 0\ar[r] & \mathfrak{X}_\infty^T \ar[r]
& \YC_\infty^T \ar[r] & \La[G]\ar[r] & \ZM_p \ar[r] & 0}}

By cutting out by any non trivial character $\psi$ of $\De$, we get an
exact sequence of $\La(\psi)[P]$-modules

\noindent\centerline{\xymatrix{ 0\ar[r] & (\mathfrak{X}_\infty^T)_\psi \ar[r]
& (\YC_\infty^T)_\psi \ar[r] & \La(\psi)[P]\ar[r] & 0}}

\noindent in which the last two modules have projective dimension less
than or equal to $1$. Then classical relations between projective
dimension in a short exact sequence (see e.g. \cite{Ng86}, proposition
3.3) show immediately the assertion $(i)$.

To go back to $(\mathfrak {X}_\infty)_\psi$, recall that Leopoldt's
weak conjecture (which is valid for the cyclotomic $\ZM_p$ extension)
yields an exact sequence

\noindent\centerline{\xymatrix{0\ar[r] & \oplus_{w\in \Si}
M_w\ar[r] & \mathfrak{X}_\infty^T \ar[r] & \mathfrak{X}_\infty^S
\ar[r] & 0}}

\noindent where $M_w$ is the module obtained by inducing $\ZM_p(1)$ from
$\Gal(F_\infty/\QM)_w$ to $\Gal(F_\infty/\QM)$
(\cite {Ng86},\cite{Wi85}). Cutting out by any locally non Teichm\"uller
character $\psi$ of $\De$ gives the assertion $(ii)$, since the
hypothesis ``locally non Teichm\"uller'' just means that $(\oplus_{w\in \Si}
M_w)_\psi =0$.
\par\qed\par
As $A:=\La(\psi)[P]$ is a local ring, lemma \ref{dpX} $(i)$ gives
us a short free resolution \xymatrix{0\ar[r]& A^m \ar[r]^-\phi &
A^m \ar[r]& (\XG_\infty^T)_\psi \ar[r] & 0 }, and $\mathrm{Fitt}_A
(\XG_\infty^T)_\psi = (\det \phi)$. Let us abbreviate $R:=\rm{Ram}
(F/\QM)$, denote $\det \phi$ by $H_R=H_{R,\psi}\in A$, and call it
``the'' {\it equivariant characteristic series} of the $A$-module
$(\XG^T_\infty)_\psi$. Lemma \ref{dpX} $(ii)$ gives us an
analogous exact sequence for $(\XG_\infty)_\psi$ and $\psi$
locally non-trivial, and analogous {\it equivariant characteristic
series}, $H=H_\psi$ of the $A$-module $(\XG_\infty)_\psi$.

Let $\psit$ be the (non-irreducible) character of $G$ induced
by $\psi$. For any $\QM^{\mathrm{sep}}_p$-irreducible character $\chi$ of
$G$ dividing $\psit$, let us denote again by $\chi$ the map obtained by
extending $\chi$ to $\ZM_p[G]$ or $\La[G]$ by linearity.
\begin{prop}\label{eqchse}
 Let $\psi$ be a non-trivial $\QM_p$-irreducible character of $\De$, then :
\begin{itemize}

\item[{\rm $(i)$}] $\mathrm{Fitt}_{\ZM_p(\psi)[P]} (\XG^T_F)_\psi=(H_R
(0))$.

\item[{\rm $(ii)$}] If moreover $\psi$ is locally non Teichm\"uller, then
$\mathrm{Fitt}_{\ZM_p(\psi)[P]} (\XG_F)_\psi=(H(0))$.

\item[{\rm $(iii)$}] For any $\QM_p^{\text {sep}}$-irreducible character
$\chi$ of $G$ dividing $\psit$~:\\
$\mathrm{Fitt}_{\La(\chi)} ((\XG_\infty^T)_\chi)=(\chi(H_R))=(h_{R,\chi}(T)),$
where $h_{R,\chi(T)}$ is the usual characteristic series of
 $(\XG_\infty^T)_\chi$ over $\La(\chi)$. Here $M_\chi$ denotes the
``$\chi$-quotient'' of $M$ (see e.g. \cite{T99},\S 2).

\item[{\rm $(iv)$}] For any $\QM_p^{\mathrm{sep}}$-irreducible character
$\chi$ of $G$~:\\
$\mathrm{Fitt}_{\La(\chi)} ((\XG_\infty)_\chi)=(\chi(H))=(h_{\chi}(T)),$
where $h_{\chi(T)}$ is the usual characteristic series of
 $(\XG_\infty)_\chi$ over $\La(\chi)$.

\end{itemize}\end{prop}
\dem  $(i)$ and $(ii)$ are direct consequences of lemma \ref{dpX}, $(iii)$ and $(iv)$ are the classical Main Conjecture (see e.g. \cite{Grei92} \cite{T99}).
\par\qed\par

Recall that by the Main Conjecture, $\chi(H) (0)\overset p \sim
h_\chi(0) \overset p \sim L_p(\chi,1)$ (see e.g. \cite{Grei92},
theorem 3.2 and proposition 3.4).

\par
 Let us now deal with
$\Sol_1(F)=\mathrm{Fitt}_{\ZM_p[G]} ((\UC_\infty/\C_\infty)_\Ga)$.
To this end, we introduce some ``structural constants'' which intervene in
the computation of Sinnott's constant ``character by character''
(see \cite{BN1})
\begin{defi}
For any $\QM_p^{\mathrm{sep}}$-irreducible character $\ch$ of $G$,
let $F_\ch$ be the field cut out by $\ch$ (i.e. the fixed field of $\Ker \ch$).
For all subfields $M\subset F$ such that $F_\ch \subset M$, define
$b_{\ch,M}^F := [F:M] \prod_{\ell} (1-\ch^{-1}(\ell))$, where the product is
taken over all primes $\ell$ dividing the conductor of $M$.
Let $b_\ch^F$ be ``the'' greatest common divisor
of all the $b_{\ch,M}^F$.\end{defi}
(actually since $O_{F_w}(\ch)$ is a local
ring, $b_{\ch}^F$ is ``equal'' to one of the $b_{\ch,M}^F$ with
minimal $p$-adic valuation). One can easily check that
the ideal $(b_\ch^F)$ defined here coincides with
$(b_\ch)$ defined in \cite{BN1}, d\'efinition 1.6.
\begin{lem}
Suppose that $p$ is totally split in $F$. For any $\ch\neq 1$, we have
 $\ch(\Sol_1(F))=(b_\ch^F L_p(\ch,1))$ in $\ZM_p(\ch)$.
\end{lem}
\dem
If $F_\ch \not \subset M$, then obviously $\ch(\sol^F_M)=0$.
By definition of $b_\ch^F$ it will be enough to show that if
$F_\ch\subset M$, then $\ch(\sol_M^F)\overset p \sim b^F_{\ch,M} L_p(\ch,1)$.
We first consider the special case $M=F$ and $\ch$ has the same
conductor, $f$ say, as $F$. In that case we have $\sol_F^F=\sol_F$ and
$b^F_{\ch,F}=1$. It follows easily that
$$
\begin{aligned}
 \ch(\sol_F) &
= \frac 1 p \sum_{g\in \Gal (F/\eQ)}
\log_p(\imath(\ep_F^g)) \ch(g^{-1}) \\&
=\frac 1 p \sum_{g\in \Gal (F/\eQ)}
\log_p(\imath(N_{\QM(\zeta_f)/F} (1-\zeta_f)^g)) \ch(g^{-1})\end{aligned}$$
$$\begin{aligned}\quad\quad &=\frac 1 p \sum_{g\in \Gal (\QM(\ze_f/\eQ))}
\log_p(\imath((1-\zeta_f)^g)) \ch(g^{-1}) \\
&=\frac 1 p \sum_{1\leq a \leq f}
\log_p(\imath(1-\zeta_f^a)) \ch^{-1}(a) \\
&\overset p \sim L_p(\ch,1),
\end{aligned}$$
(the last $p$-adic equivalence comes from a classical formula
for $L_p(\ch,1)$; see e.g \cite{Wa}, p. 63). Note that an implicit
choice of $\imath$ is made there, and that since we are assuming
that $p$ is not ramified in $F$, the quantities $p-\ch(p)$, $f_\ch$, and
$\tau(\ch)$ are all $p$-adic units.
\par
Now $\ch (\widetilde{\sol}^F_{F_\ch})$ is well defined
because taking another lift  would only add an element
of the (additive) kernel of $\ch$.
From the previous special case we deduce the equivalence (without
assuming anything on the conductor $f_\ch$) :
\noindent $\ch (\widetilde{\sol}^F_{F_\ch})=
\ch(\sol^{F_\ch}_{F_\ch})\overset p \sim L_p(\ch,1)$, hence
$\ch(\sol^F_{F_\ch})\sim \ch(\widetilde{\sol}^F_{F_\ch}\ {\Tr}_{F/F_\ch})
\sim [F:F_\ch] \ch(\widetilde{\sol}^F_{F_\ch})
  \sim [F:F_\chi] L_p(\ch,1) \sim b^F_{\ch,F_\ch} L_p(\ch,1).$
The remaining cases are treated using the distribution relations satisfied
by the cyclotomic numbers (and hence by the elements $\sol_M^F$)
and the formal identity $\ch(\sol^F_M)=\frac 1 {[M:F_\ch]}
 \ch({\Tr}_{M/F_\ch} \sol^F_M)$.\par\qed
\begin{cor}\label{compare} $MW(F)\neq \Sol_1(F)$ in general. But in the
semi-simple case (i.e. $p\nmid [F:\QM]$), $MW(F)=\Sol_1(F)$.\end{cor}
\dem In the semi-simple case, for any $\ch \neq 1$, $b_\ch^F$ is a $p$-adic
 unit, the preceding calculation about $P$-cohomology becomes empty.
As for the trivial character, the corresponding idempotent is just the
norm map, which brings us down to $\QM$, where nothing harmful happens.
\par\qed
\begin{theo} Suppose that $p$ is totally split in $F$. For any non
trivial and locally
non Teichm\"uller $\QM_p$-irreducible character $\psi$ of $\De$,
the following properties
are equivalent :
\begin{itemize}
\item[{\rm $(i)$}] $(\C_\infty)_\psi$ is $\La(\psi)[P]$-free
\item[{\rm $(ii)$}] $(\XG_\infty)_\psi$ and $(\UC_\infty/\C_\infty)_\psi$ have
``the'' same equivariant characteristic series in $\La (\psi) [P]$
\item[{\rm $(iii)$}] $(KN_F)_\psi=0$ and
$(\pucod_F)_\psi $ is $\ZM_p(\psi) [P]$-free
\end{itemize}
They all imply :
\begin{itemize}
\item[{\rm $(iv)$}] $MW(F)_\psi=\Sol_1(F)_\psi$
\item[{\rm $(v)$}] $(\Uha^\prime_F/\pucod_F)_\psi^\ast$ and
$(X^\prime_\infty)_{\Ga,\psi}$ have the same Fitting ideal over
$\ZM_p(\psi)[P]$, $(.)^\ast$ denoting the Pontryagin dual, endowed
with the Galois action defined by $f^\si(x)=f(\si x)$.
\end{itemize}\end{theo}
\dem To show the equivalence between $(i)$ and $(ii)$, let us consider
the module $(\UC_\infty/\C_\infty)_\psi \simeq (\UC_\infty)_\psi/
(\C_\infty)_\psi$. Since $p$ is totally split in $F$, for all $v\in S$ the
local fields $F_v$ contain no $p^{\text{th}}$ power root of unity.
Coleman's theory (see e.g. the exact sequence in
theorem 4.2 of \cite{T99}) then shows that $\UC_\infty$ is a rank one
$\La[G]$-free module, hence $(\UC_\infty)_\psi
\simeq A:=\La(\psi)[P]$. This gives the short resolution
\xymatrix{
0\ar[r] & (\C_\infty)_\psi \ar[r] & A
\ar[r] & (\UC_\infty/\C_\infty)_\psi\ar[r]
& 0 }, which shows the equivalences~:
\par
\noindent $\mathrm{pd}_R (\UC_\infty/\C_\infty)_\psi \leq 1
\iff (\C_\infty)_\psi$ is $R$-free $\iff \mathrm{Fitt}_R
(\UC_\infty/\C_\infty)_\psi$ is principal.
\par
In the last eventuality, let us denote by $J$ the equivariant characteristic
series of $(\UC_\infty/\C_\infty)_\psi$ (as defined before \ref{eqchse}).
To compare $J$ with $H$ (the equivariant characteristic series of
$(\XG_\infty)_\psi$), we appeal to still another algebraic lemma of Greither~:
\begin{lem}[\cite{Grei98}, 3.7] Let $A=\La(\psi) [P]$. If $M$ is an $A$-torsion
module of projective dimension at most $1$, such that $M/pM$ is
finite, and if $J$ is an element of $A$ such that
$\mathrm{Fitt}_{\La(\ch)} (M_\ch) = (\ch(J))$ for all
$\QM_p^{\mathrm{sep}}$-irreducible characters $\ch$ of $G$
dividing $\widetilde\psi$, then actually $\mathrm{Fitt}_A(M)=(J)$.
\end{lem}
We apply this to $M=(\UC_\infty/\C_\infty)_\psi$. It is known (see \cite{BN1},
 \cite{T99}) that for all $\ch\mid\psit$, the $\La(\ch)$-modules
$(\UC_\infty/\C_\infty)_\ch$ and $(\XG_\infty)_\ch$ have the same (usual)
characteristic series. In particular, their $\mu$-invariants are null.
It follows that $\ch(H)=\ch(J)$ for all $\ch \mid \psit$, hence
$(H)=(J)$ by Greither's lemma.
Conversely the equality $(H)=\mathrm{Fitt}_{\La(\psi)[P]}
(\UC_\infty/\C_\infty)_\psi$
implies the principality of this last ideal.
The proof of the equivalence between $(i)$ and $(ii)$ is thus complete.
Property
$(i)$ implies the triviality of $(KN_F)_\psi$ and then $(iii)$
by $\Ga$-co-descent.
\par
Conversely, assume $(iii)$ and choose a $\ZM_p(\psi)[P]$-basis
$(\overline {y_i})$ of $(\pucod_F)_\psi$. By Nakayama's lemma, this can be
lifted to a system of $\La(\psi)[P]$-generators $(y_i)$ of
$(\C^{\prime\prime}_\infty)_\psi$. Any linear relation $\sum \la_i y_i=0$,
with $\la_i \in \La(\psi)[P]$, would give, by $\Ga$-co-descent,
$\overline{\la_i}=0$ in $\ZM_p(\psi)[P]$ for any $i$, viz. $T(:=\ga-1)$
would divide all the coefficients $\la_i$. As $(\C_\infty)_\psi$ has
no $T$-torsion, we could then simplify by $T$ in the above linear
relation and repeat the process. This shows that the system $(y_i)$
is a $\La(\psi)[P]$-basis.
\par
By the same argument, $(ii)$ implies $(iv)$.
Property $(ii)$ implies $(v)$ because of the $(\psi)$-parts of the exact
sequence $(7)$ together with one last algebraic lemma~:
\begin{lem}\label{CGle}
Let {\xymatrix@=11.5pt{0\ar[r] & M \ar[r] & N \ar[r] & N^\prime
\ar[r] & M^\prime \ar[r] & 0}} be an exact sequence of
$\ZM_p(\psi)[P]$-modules of finite order. Suppose that $N$ and
$N^\prime$ are of projective dimension at most $1$ and have the
same Fitting ideal over $\ZM_p(\psi)[P]$. Then $M^\ast$ and $M^\prime$
have the same Fitting ideal over $\ZM_p(\psi)[P]$.
\end{lem}
\dem This is \cite{CG98}, Proposition 6, but note that the correct
statement involves $M^\ast$ (and not $M$ as in \cite{CG98}).
\par\qed\par

{\bf Acknowledgement} : We thank Radan Ku\v cera for pointing out an error
in a previous formulation of lemma \ref{Rub}.

\par
\centerline{ \hfill {Jean-Robert Belliard} \hfill \hfill {Th{\cfac{o}}ng
Nguy{\cftil{e}}n-Quang-{\Dbar}{\cftil{o}}} \hfill }
\par
\centerline { \hfill {belliard@math.univ-fcomte.fr} \hfill\hfill
{nguyen@math.univ-fcomte.fr} \hfill }
\par
\centerline{Universit\'e de Franche-Comt\'e,}
\par
\centerline{UMR 6623,}
\par
\centerline{16, route de Gray,}
\par
\centerline{25030 Besan\c con cedex.}
\par
\centerline{FRANCE}
\end{document}

%% file: engladef.tex
\newtheorem{theo}{Theorem}[section]
\newtheorem{prop}[theo]{Proposition}
\newtheorem{cor}[theo]{Corollary}
\newtheorem{lem}[theo]{Lemma}
\newtheorem{defi}[theo]{Definition}
\newtheorem{con}[theo]{Conjecture}

\newtheorem{nots}[theo]{Notations}
\newtheorem{nota}[theo]{Notation}

\def \fr {\mathrm {fr}}
\def \demdu#1 { {\sl Proof #1.} }

\def \Gal {\text{\rm Gal}}
\def \Cal {\mathcal}

\def \tor {\text {\rm Tor}}
\def \coker {\Coker}
\def \dem{{\sl Proof.} }
\def \eZ{{\mathbb {Z}} }
\def \eQ{{\mathbb {Q}} }
\def \limpro{\lim\limits_{\leftarrow} }

\def \qed{\text {$\quad \square$}}

\def \C {\overline {C}}
\def \U {\overline {U}}
\def \Dnb {\text {\rm Cyc}}
\def \Dnbp {\text{\rm Cyc}^\prime}

\def \F#1 {\overline {F^\times_{#1}}}

\def \Cl {\overline{\text{Cl}}}

\def \pucod {\widetilde {C}^{\prime\prime}}

\def \rle {R^{\text{Leop}}}
\def \rgr {R^{\text{Gross}}}
\def \D {\overline {\Dnb}}
\def \Dp {\overline {\Dnbp}}

\def \cond {\mathrm {cond}}

\def\BM{{\mathbb{B}}}
\def\CM{{\mathbb{C}}}

\def\NM{{\mathbb{N}}}

\def\QM{{\mathbb{Q}}}

\def\ZM{{\mathbb{Z}}}

\def\LG{{\mathfrak L}}

\def\PG{{\mathfrak P}}

\def\XG{{\mathfrak X}}

\def\al{\alpha}
\def\be{\beta}
\def\ga{\gamma}
\def\Ga{\Gamma}
\def\de{\delta}
\def\De{\Delta}
\def\ep{\varepsilon}

\def\ch{\chi}
\def\la{\lambda}
\def\La{\Lambda}

\def\rh{\rho}
\def\si{\sigma}
\def\Si{\Sigma}
\def\th{\theta}
\def\Th{\Theta}

\def\ze{\zeta}

\def\CC{{\mathcal{C}}}
\def\DC{{\mathcal{D}}}

\def\FC{{\mathcal{F}}}

\def\KC{{\mathcal{K}}}

\def\OC{{\mathcal{O}}}

\def\SC{{\mathcal{S}}}

\def\UC{{\mathcal{U}}}

\def\YC{{\mathcal{Y}}}

\def\Lti{{\widetilde{L}}}

\def\Fha{{\widehat{F}}}

\def\Uha{{\widehat{U}}}

\def\Fba{{\overline{F}}}

\def\psit{{\widetilde{\psi}}}

\def\stc{\mathop{\mathrm{Stick}}\nolimits}

\def\Cl{\mathop{\mathrm{Cl}}\nolimits}

\def\Coker{\mathop{\mathrm{Coker}}\nolimits}

\def\Gal{\mathop{\mathrm{Gal}}\nolimits}

\def\Im{\mathop{\mathrm{Im}}\nolimits}
\def\im{\mathop{\mathrm{Im}}\nolimits}

\def\Ker{\mathop{\mathrm{Ker}}\nolimits}

\def\log{\mathop{\mathrm{log}}\nolimits}

\def\Tr{\mathop{\mathrm{Tr}}\nolimits}

\def\pr{\prime}